\documentclass{eajam}
\renewcommand\thisnumber{x}
\renewcommand\thisyear {200x}
\renewcommand\thismonth{xxx}
\renewcommand\thisvolume{xx}

\usepackage{charter}
\usepackage[charter]{mathdesign}

\usepackage{subfig}
\usepackage{booktabs}
\usepackage{algorithm}
\usepackage{algpseudocode}
\usepackage{xcolor}

\numberwithin{equation}{section}

\begin{document}

\markboth{Rui Sheng, Peiying Wu, Jerry Zhijian Yang and Cheng Yuan}{Solving the inverse source problem by MC-fPINNs}
\title{Solving the inverse source problem of the fractional Poisson equation by MC-fPINNs}


\author[Authors]{Rui Sheng\affil{1}, Peiying Wu\affil{1}, Jerry Zhijian Yang\affil{1,2} and Cheng Yuan\affil{3,4}\comma\corrauth}
\address{\affilnum{1}\ School of Mathematics and Statistics, Wuhan University, Wuhan, 430072, China\\
\affilnum{2}\ Hubei Key Laboratory of Computational Science, Wuhan University, Wuhan, 430072, China\\
\affilnum{3}\ School of Mathematics and Statistics, and Hubei Key Lab--Math. Sci., Central China Normal University, Wuhan 430079, China\\
\affilnum{4}\ Key Laboratory of Nonlinear Analysis $\&$ Applications (Ministry of Education), Central China Normal University, Wuhan 430079, China\\}

%
\emails{{\tt shengrui@whu.edu.cn} (R. Sheng), {\tt peiyingwu@whu.edu.cn} (P. Wu), {\tt zjyang.math@whu.edu.cn} (J. Z. Yang), {\tt yuancheng@ccnu.edu.cn} (C. Yuan)}
%
\begin{abstract}
In this paper, we effectively solve the inverse source problem of the fractional Poisson equation using MC-fPINNs. We construct two neural networks $ u_{NN}(x;\theta )$ and $f_{NN}(x;\psi)$ to approximate the solution $u^{*}(x)$ and the forcing term $f^{*}(x)$ of the fractional Poisson equation. To optimize these two neural networks, we use the Monte Carlo sampling method mentioned in MC-fPINNs and define a new loss function combining measurement data and the underlying physical model. Meanwhile, we present a comprehensive error analysis for this method, along with a prior rule to select the appropriate parameters of neural networks. Several numerical examples are given to demonstrate the great precision and robustness of this method in solving high-dimensional problems up to 10D, with various fractional order $\alpha$ and different noise levels of the measurement data ranging from 1$\%$ to 10$\%$. 
\end{abstract}



\maketitle


\section{Introduction\label{sec1}}
The fractional partial differential equations (fPDEs) has played a critical role in the mathematical modeling of anomalous phenomena in fields of science and engineering, such as hydrology~\cite{uchaikin2013fractional}, viscoelasticity \cite{mainardi2022fractional} and turbulent flow \cite{2006turbulence, epps2018turbulence}. In practical applications, we often need to recover some information missed in these fPDEs, including the diffusion coefficients, initial or boundary data, or source terms, especially for some problems in physics \cite{FAN201540} and biology \cite{LIU2011822}. Among these, due to the non-local property and singularity of the fractional Laplacian operators, solving the inverse problem of the fractional Poisson equation is still a difficult task.

More recently, with the significant advancements in deep learning techniques applied to computational vision and natural language processing, neural network architectures have also been employed for solving the inverse problem of fPDEs. For example, Pang et al. extended PINNs~\cite{raissi2019physics} to fPINNs~\cite{pang2019fpinns} to solve the inverse problem of identifying parameters in the partial advection diffusion equation. In their method, the integer derivative is calculated using automatic differentiation while the fractional derivative is approximated by traditional numerical discretization, which leads to a high computational cost especially for high-dimensional problems. Later Yan et al. \cite{Yan2022BayesianIW} proposed a Bayesian Inversion with Neural Operator (BINO) approach to solve the Bayesian inverse problem of the time fractional subdiffusion equation, of which the diffusion coefficient can be recovered. On the other hand, they \cite{yan2023laplacefpinns} presented the Laplace fPINNs method to identify the diffusion coefficient function in time-fractional diffusion equation. This method first transforms the original equation into a restricted problem in Laplace space and then solves it using PINNs. However, although BINO and Laplace-fPINNs have been used successfully to solve the time-fractional problems, they did not consider the space-fractional and high-dimensional problems.

To deal with space-fractional derivative and higher-dimensional problems more efficiently, Guo et al. \cite{GUO2022115523} propose a Monte Carlo sampling-based PINNs method (MC-fPINN) to identify the parameters in the fractional advection-diffusion equation. Different from fPINNs, MC-fPINNs calculate fractional derivatives using Monte Carlo sampling instead of traditional methods, resulting in reduced computational cost and enabling their application to higher-dimensional problems. In this paper, we will extend the MC-fPINNs method to address the inverse source problem of the fractional Poisson equation, not only in low dimensions but also in high dimensions. In summary, our main contributions are as follows.  

\begin{itemize}
\item 
Except for the neural network for approximating the solutions of the fractional Poisson equation, we represent the forcing term function as another fully-connected neural network. To optimize these two neural networks, we use the Monte Carlo sampling method mentioned in MC-fPINNs and define a new loss function containing regularization terms for residuals of fPDEs and measurement data. We effectively address the inverse source problems related to the fractional Poisson equation, even in high-dimensional cases. In our computational experiments, we show the high accuracy and reliability of our method by varying the fractional order $\alpha$ while introducing noise to the measurement data at levels between 1$\%$ and 10$\%$.
\item
We provide a systemetical error analysis for the MC-fPINNs method for this inverse problem, as well as a guideline for choosing suitable neural network parameters, such as the total number of non-zero weights, depth and samples. This result guarantees that the error in estimation remains consistent and manageable in order to attain the desired level of accuracy.
\end{itemize}

The rest of this paper is organized as follows. In Section \ref{Pre}, we describe the inverse problem of the fractional Poisson equation with Dirichlet boundary condition, and the notations that will be used in this work. In Section \ref{Method}, the Inverse MC-fPINNs method will be introduced. In Section \ref{error}, the error analysis and convergence rate of MC-fPINNs for the inverse source problem are given. In Section \ref{experiment}, we provide some numerical simulations. A conclusion and short discussion is given in Section \ref{conclusion}.

\section{Preliminary}\label{Pre}
\subsection{Problem Setup}\label{Setup}
In this paper, we consider the inverse problem for the following fractional Poisson equation with homogeneous Dirichlet boundary condition:
\begin{equation}\label{eq:problemsetup}
\begin{aligned}
(-\Delta)^{\alpha/2}u(x) ={f}(x),  \qquad  &x\in \Omega,\\
u(x)=0,  \qquad \qquad  &x\in R^d \backslash\Omega,
\end{aligned}
\end{equation}
where $\Omega \subseteq R^{d}$  is a bounded domain with smooth boundary $\partial \Omega$, and the fractional laplacian operator is defined by:
\begin{equation}\label{def:flap}
(-\Delta)^{\alpha/2}u(x)\triangleq C(d,\alpha)P.V.\int_{R^{d}} \frac{u(x)-u(y)}{\left \| x-y \right \|_2^{d+\alpha}} dy,\quad 0<\alpha<2.
\end{equation}
Here "P.V."denotes the principal value of the integral and the constant $C(d,\alpha)$ is given by
\begin{equation}\label{definition:C(d,s)}
C(d,\alpha )=\frac{2^\alpha \Gamma(\frac{\alpha+d}{2})}{\pi^{d/2}\left | \Gamma(-\alpha/2) \right | }.
\end{equation}
Without loss of generality, we may assume that $ \Omega \subseteq {\left [0,1\right ]^d}$. In this work, our objective is to simultaneously recover both the forcing term $f(x)$ and the solution $u(x)$, given the boundary condition and some extra measurements of $u$.

\subsection{Notations of Neural Networks and Fractional Function Spaces}
A $\varrho$-neural network $\phi:\mathbb{R}^{N_{0}}\rightarrow\mathbb{R}^{N_{\mathcal{D}+1}}$ refers to a function defined by 
\begin{equation*}
\phi(x)=T_{\mathcal{D}}(\varrho(T_{\mathcal{D}-1}(\cdots\varrho(T_{0}(x))\cdots))),
\end{equation*}
where the activation function $\varrho$ is applied component-wisely and $T_{\ell}(x):=A_{\ell}x+b_{\ell}$ is an affine transformation with $A_{\ell}\in\mathbb{R}^{N_{\ell+1}\times N_{\ell}}$ and $b_{\ell}\in\mathbb{R}^{N_{\ell}}$ for $\ell=0,\ldots\mathcal{D}$. In this paper, we consider the cases $N_{0}=d$ and $N_{\mathcal{D}+1}=1$. The number $\mathcal{D}$ the depth of the neural network. We denote $\mathcal{S}_{i}=\sum_{\ell=1}^{i}(\|A_{\ell}\|_{0}+\|b_{\ell}\|_{0})$ for $i=1,\cdots,\mathcal{D}$ as the number of non-zero weights in the first $i$-layers and $\mathcal{S}=\mathcal{S}_{\mathcal{D}}$ as the total number of nonzero weights. Based on these, we may further define$\mathcal{N}_{\varrho}(\mathcal{D},\mathcal{S},\mathcal{B})$ as the collection of $\varrho$-neural networks with at most $\mathcal{D}$ layers and $\mathcal{S}$ non-zero weights, and each weight is bounded by $\mathcal{B}$.

Next, according to \cite{antil2017fractional}\cite{dinezza2011hitchhikers}\cite{ROSOTON2014275}, the fractional function spaces with the corresponding norms used in this paper will be introduced as follows. Let $\alpha=\left(\alpha_{1}, \cdots, \alpha_{n}\right) $ be an $n-$dimensional index vector with $|\alpha|:=\sum_{i=1}^{n} \alpha_{i}$. For any real $s>0$ and $p \in  [1,+\infty)$, by setting $\theta=s-\left \lfloor s  \right \rfloor $ with $m=\left \lfloor s  \right \rfloor$, we can define the Sobolev space
\begin{equation*}
W^{s, p}(\Omega):=\left\{u \in W^{m, p}(\Omega)\left|\int_{\Omega} \int_{\Omega} \frac{\left|D^{\alpha} u(x)-D^{\alpha} u(y)\right|^{p}}{|x-y|^{\theta p+d}} d x d y<\infty,   \quad \forall |\alpha|=\lfloor s\rfloor\right.\right\},
\end{equation*}

\begin{equation*}
\|u\|_{W^{s, p}(\Omega)}=\left(\|u\|_{W^{m, p}(\Omega)}^{p}+\sum_{|\alpha|=m} \int_{\Omega} \int_{\Omega} \frac{\left|D^{\alpha} u(x)-D^{\alpha} u(y)\right|^{p}}{|x-y|^{\theta p+d}} d x d y\right)^{1 / p},
\end{equation*}
while the term 
\begin{equation*}
[u]_{W^{s, p}(\Omega)}:=\left(\int_{\Omega} \int_{\Omega} \frac{|u(x)-u(y)|^{p}}{|x-y|^{n+s p}} d x d y\right)^{\frac{1}{p}}
\end{equation*}
is the Gagliardo (semi)norm of u. When $p=2$, $W^{s, p}(\Omega)$ is a fractional Hilbert Space which can be denoted by $H^s(\Omega)$. And we use the notation $C^s(\Omega)$ to refer to the Hölder space $C^{m,\theta}$:
\begin{equation*}
C^s(\Omega)=C^{m, \theta}(\Omega):=\left\{u \in C^{m}(\Omega) \mid\|u\|_{C^{m, \theta}(\Omega)}<\infty\right\},
\end{equation*}
where the norm $\|\cdot\|_{C^{m, \theta}(\Omega)}$ is defined as
\begin{equation*}
\|u\|_{C^{s}(\Omega)}=\|u\|_{C^{m, \theta}(\Omega)}=\sum_{|\alpha| \leq m}\left\|D^{\alpha} u\right\|_{C^{0}(\Omega)}+\sum_{|\alpha|=m}\sup _{x, y \in \Omega, x \neq y} \frac{|D^{\alpha}u(x)-D^{\alpha}u(y)|}{|x-y|^{\theta}} .
\end{equation*}

\section{Methodology}\label{Method}
\subsection{Inverse MC-fPINNs Method}
To solve the inverse problem introduced in subsection \ref{Setup}, we define the following loss function:
\begin{equation}\label{eq:loss}
\mathcal{L}(u,f)= \left \| L(u)-f \right \|_{L^2(\Omega)}^2+\left \| u\right \|_{L^2(\partial \Omega)}^2+\left \| u-u_d \right \|_{L^2(\Omega)}^2,
\end{equation}
where $L(u):=(-\Delta)^{\alpha/2}u$ and $u_d$ is the extra measurement of u. To reconstruct the solution $u$ and the force term $f$ simultaneously, we would use $u_{NN}(x,\theta)\in\mathcal{N}_u:=\mathcal{N}_{\rho}\left(\mathcal{D}_u, \mathfrak{n}_{\mathcal{D}_u}, B_{\theta}\right)$  and  $f_{NN}(x;\psi)\in \mathcal{N}_f:=\mathcal{N}_{\rho}\left(\mathcal{D}_f, \mathfrak{n}_{\mathcal{D}_f}, B_{\psi}\right)$ to approximate them respectively. Here we adopt $tanh(x)$ as the activation function $\rho$.

For the fractional Laplacian, the integral in \eqref{def:flap} can be divided as in MC-fPINNs \cite{GUO2022115523}:
\begin{equation}\label{equ:laplace decomposition}
\begin{aligned}
L(u_{N N}(x))&=(-\Delta)^{\alpha / 2} u_{N N}(x)\\
&=C_{d, \alpha}\left[\underbrace{\int_{y \in B_{r_{0}}(x)} \frac{u_{N N}(x)-u_{N N}(y)}{\|x-y\|_{2}^{d+\alpha}} \mathrm{d} y}_{I}+\underbrace{\int_{y \notin B_{r_{0}}(x)} \frac{u_{N N}(x)-u_{N N}(y)}{\|x-y\|_{2}^{d+\alpha}} \mathrm{d} y}_{II}\right].
\end{aligned}
\end{equation}
where $B_{r_{0}}(x):=\left\{y \mid\|y-x\| \leq r_{0}\right\}$. The first part can be calculated as
\begin{equation}\label{equ:decomposition:1}
\begin{aligned}
I&=\frac{\left|S^{d-1}\right| r_{0}^{2-\alpha}}{2(2-\alpha)} \mathbb{E}_{\xi, r_{I} \sim f_{I}\left(r_{I}\right)}\left[\frac{2 u_{N N}(x)-u_{N N}\left(x-r_{I} \xi\right)-u_{N N}\left(x+r_{I} \xi\right)}{r_{I}^{2}}\right],\\
&\approx \frac{\left|S^{d-1}\right| r_{0}^{2-\alpha}}{2(2-\alpha)} \mathbb{E}_{\xi, r_{I} \sim f_{I}\left(r_{I}\right)}\left[\frac{2 u_{N N}(x)-u_{N N}\left(x-r_{\epsilon} \xi\right)-u_{N N}\left(x+r_{\epsilon} \xi\right)}{r_{\epsilon}^{2}}\right],
\end{aligned}
\end{equation}
where $\xi$  is a random variable uniformly distributed on the unit $(d-1)$-sphere $S^{d-1}$, $\left|S^{d-1}\right|$ denotes the surface area of  $S^{d-1}$,
$f_{I}\left(r_{I}\right)=\frac{2-\alpha}{r_{0}^{2-\alpha}} r_{I}^{1-\alpha} \cdot 1_{r_{I} \in\left[0, r_{0}\right]}$ and $r_\epsilon=\max\left\{\epsilon,r_I\right \}$. $r_I$ can be sampled as:
\begin{equation*}
r_I/r_0\sim \text{Beta}(2-\alpha,1).
\end{equation*}
As for the second part,we have:
\begin{equation}\label{equ:decomposition:2}
II=\frac{\left|S^{d-1}\right| r_{0}^{-\alpha}}{2 \alpha} \mathbb{E}_{\xi, r_{o} \sim f_{O}\left(r_{o}\right)}\left[2 u_{N N}(x)-u_{N N}\left(x-r_{o} \xi\right)-u_{N N}\left(x+r_{o} \xi\right)\right]
\end{equation}
where $f_{O}\left(r_{o}\right)=\alpha r_{0}^{\alpha} r_{o}^{-1-\alpha} 1_{r_{o} \in\left[r_{0}, \infty\right)}$ and $r_o$ can be sampled as:
\begin{equation*}
r_0/r_o\sim \text{Beta}(\alpha,1).
\end{equation*}

Now we focus on the discretization of $\mathcal{L}(u,f)$. With the Monte Carlo method, the above integral can be approximated by
\begin{eqnarray*}
  \mathit{L}(u_{NN}(x;\theta)) &\approx & \hat{\mathit{L}}(u_{NN}(x;\theta);r_\epsilon,r_o,\xi), \\
 & = & \frac{1}{m} \sum_{j=1}^{m} (\frac{C_{d,\alpha }\left | S^{d-1} \right |{r_{0}}^{2-\alpha} }{2(2-\alpha)}\cdot \frac{2u_{NN}(x;\theta)
-u_{NN}(x-r_{\epsilon j} \xi_j;\theta)-u_{NN}(x+r_{\epsilon j} \xi_j;\theta)}{r_{\epsilon j}^2} \\
&&+\frac{C_{d,\alpha }\left | S^{d-1} \right |{r_0}^{-\alpha}}{2\alpha}\cdot (2u_{NN}(x;\theta)-u_{NN}(x-r_{oj} \xi_j;\theta)-u_{NN}(x+r_{oj} \xi_j;\theta)),
\end{eqnarray*}
in which $\left \{ r_{\epsilon j}, r_{oj},\xi_j \right \}_{j=1}^m$ is a group of parameters sampled according to their corresponding distributions. Furthermore, to achieve an unbiased estimation of the residual $\| L(u)-f \|_{L^2(\Omega)}^2$, we would sample another group of parameters $\left \{ {r_{\epsilon j}}', {r_{oj}}', {\xi_j}' \right \}_{j=1}^m$  independently according to the same distributions. Let 
\begin{equation*}
  P_j=\left \{ r_{\epsilon j}, r_{oj}, \xi_j, {r_{\epsilon j}}', {r_{oj}}',{\xi_j}' \right \}, \quad 1\leq j\leq m,
\end{equation*}
we may define the unbiased estimation as
\begin{equation*}
\mathcal{L}{equ}(\theta,\psi) =|\Omega|\frac{1}{N}\frac{1}{m}\sum_{i=1}^{N}\sum_{j=1}^{m} \mu(x_i,P_j)\cdot\eta(x_i,P_j),
\end{equation*}
where
\begin{equation*}
\begin{aligned}
\mu(x,P)=&\frac{C_{d,\alpha }\left | S^{d-1} \right |{r_{0}}^{2-\alpha} }{2(2-\alpha)}\cdot \frac{2u_{NN}(x;\theta)
-u_{NN}(x-r_{\epsilon} \xi;\theta)-u_{NN}(x+r_{\epsilon } \xi;\theta)}{r_{\epsilon }^2}\\
&+\frac{C_{d,\alpha }\left | S^{d-1} \right |{r_0}^{-\alpha}}{2\alpha}\cdot (2u_{NN}(x;\theta)-u_{NN}(x-r_{o} \xi;\theta)-u_{NN}(x+r_{o} \xi;\theta)-f_{NN}(x)),
\end{aligned}
\end{equation*}
\begin{equation*}
\begin{aligned}
\eta(x,P)=&\frac{C_{d,\alpha }\left | S^{d-1} \right |{r_{0}}^{2-\alpha} }{2(2-\alpha)}\cdot \frac{2u_{NN}(x;\theta)
-u_{NN}(x-{r_{\epsilon j}}' {\xi}';\theta)-u_{NN}(x+{r_{\epsilon }}' {\xi}';\theta)}{{r_{\epsilon }}'^2}\\
&+\frac{C_{d,\alpha }\left | S^{d-1} \right |{r_0}^{-\alpha}}{2\alpha}\cdot (2u_{NN}(x;\theta)-u_{NN}(x-{r_{o}}' {\xi}' ;\theta)-u_{NN}(x+{r_{o}}' {\xi}';\theta)-f_{NN}(x)).
\end{aligned}
\end{equation*}
On the other hand, the remaining two terms in loss can also be approximated by
\begin{equation*}
\mathcal{L}_g(\theta) = |\partial \Omega|\frac{1}{N_g} \sum_{i=1}^{N_g} {[u_{NN}(x_i,\theta))-0]}^2, 
\end{equation*}
\begin{equation*}
\mathcal{L}_u(\theta) = |\Omega|\frac{1}{N_u} \sum_{i=1}^{N_u} {[ u_{NN}(x_i,\theta))-u_i]}^2, 
\end{equation*}
with which we may obtain the following empirical risk
\begin{equation}\label{eq:empirical loss}
\hat{\mathcal{L}}(\theta,\psi)=\omega_{equ}\mathcal{L}{equ}(\theta,\psi)+\omega_g \mathcal{L}_g+\omega_u \mathcal{L}_u.
\end{equation}

Here $\omega_{equ}$, $\omega_g$ and $\omega_u$ represent the weights of different loss terms. With the help of some optimization algorithm, the parameters $\psi$ and $\theta$ can be learned by minimizing the loss. The complete process is illustrated in Algorithm \ref{alg:MC}.

\begin{algorithm}[H]
\caption{MC-fPINNs for Inverse Source Problem of the fractional Poisson Equation.}
\label{alg:MC}
\begin{algorithmic}[1]
\Require Prepare training data $\mathcal{D}=\{(x_{i})\}_{i=1}^{N_{equ}}\cup\{{x}_{i}\}_{i=1}^{N_g}$ and measurement data set $S=\{(x_{i},u_{i})\}_{i=1}^{N_u}$.
\State Construct neural networks $(f_{\psi},u_{\theta})$, parameterized by $(\psi,\theta)$.
\State Initialize parameters $(\psi,\theta)$ randomly.
\For {$epoch=1:K$}
\State {$\hat{\mathcal{L}}_n(f_{\psi},u_{\theta})=\omega_{equ}\mathcal{L}{equ}(f_{\psi},u_{\theta})+\omega_g \mathcal{L}_g(u_{\theta})+\omega_u \mathcal{L}_u(u_{\theta})$. as (\ref{eq:empirical loss})}
\State {$(g_{\psi},g_{\theta})=\nabla_{(\psi,\theta)}{\hat{\mathcal{L}}_n}(f_{\psi},u_{\theta})$.}
 \State {$(\psi,\theta)\leftarrow\operatorname{SGD}\big\{(\psi,\theta),(g_{\psi},g_{\theta}),\tau\big\}$.}
\EndFor
\Ensure {$\tilde{f}_{n}=f_{\psi}$, $\tilde{u}_{n}=u_{\theta}$.}
\end{algorithmic}
\end{algorithm}
\section{Error Analysis}\label{error}
\subsection{Error Decomposition}
In the following we decompose the error into the approximation error and statistical error.
\begin{proposition}\label{proposition:error-decomposition}
Let$(u^*,f^*)\in\underset{u\in\mathcal{U},f\in\mathcal{F}}{argmin} \mathcal{L}(u,f)$ and $(\hat{u},\hat{f)}\in \underset{u\in\mathcal{N}_u,f\in\mathcal{N}_f}{argmin}{\hat{\mathcal{L}}(u,f)}$, $\mathcal{L}$ and $\hat{\mathcal{L}}$ are defined in (\ref{eq:loss}) and (\ref{eq:empirical loss}) respectively. Then for $\forall\bar{u}\in\mathcal{N}_u,\forall\bar{f}\in\mathcal{N}_f$,we have:
\begin{equation*}
\mathcal{L}(\hat{u},\hat{f})-\mathcal{L}(u^*,f^*)\le \underset{\varepsilon_{app}}{\underbrace{\underset{\bar{u}\in\mathcal{N}_u ,\bar{f}\in\mathcal{N}_f}{\inf}[\mathcal{L}(\bar{u},\bar{f})-\mathcal{L}(u^*,f^*)]}}+\underset{ \varepsilon_{sta} }{\underbrace{2\underset{u\in\mathcal{N}_u ,f\in\mathcal{N}_f}{\sup}\left |\hat{\mathcal{L}}(u,f)-\mathcal{L}(u,f)\right |}}. 
\end{equation*}
\end{proposition}
\begin{proof}
For $\forall\bar{u}\in\mathcal{N}_u,\forall\bar{f}\in\mathcal{N}_f$,we obtain:
\begin{equation*}
\begin{aligned}
\mathcal{L}(\hat{u},\hat{f})-\mathcal{L}(u^*,f^*)=&[\mathcal{L}(\bar{u},\bar{f})-\mathcal{L}(u^*,f^*)]+[\hat{\mathcal{L}}(\bar{u},\bar{f})-(\mathcal{L}(\bar{u},\bar{f})]\\
&-[\hat{\mathcal{L}}(\bar{u},\bar{f})-\hat{\mathcal{L}}(\hat{u},\hat{f})]+[\mathcal{L}(\hat{u},\hat{f})]-\hat{\mathcal{L}}(\hat{u},\hat{f})].
\end{aligned}
\end{equation*}
Since $\hat{\mathcal{L}}(\bar{u},\bar{f})-\hat{\mathcal{L}}(\hat{u},\hat{f})\ge 0$, then
\begin{equation*}
\begin{aligned}
\mathcal{L}(\hat{u},\hat{f})-\mathcal{L}(u^*,f^*)&\le [\mathcal{L}(\bar{u},\bar{f})-\mathcal{L}(u^*,f^*)]+\left | \mathcal{L}(\bar{u},\bar{f})-\hat{\mathcal{L}}(\bar{u},\bar{f}) \right | 
+\left |\hat{\mathcal{L}}(\hat{u},\hat{f})-\mathcal{L}(\hat{u},\hat{f})\right | \\
&\le [\mathcal{L}(\bar{u},\bar{f})-\mathcal{L}(u^*,f^*)]+2\underset{u\in\mathcal{N}_u ,f\in\mathcal{N}_f}{\sup}\left |\hat{\mathcal{L}}(u,f)-\mathcal{L}(u,f)\right |.
\end{aligned}
\end{equation*}
Next we perform $\underset{\bar{u}\in\mathcal{N}_u ,\bar{f}\in\mathcal{N}_f}{\inf}$ on both sides, then 
\begin{equation*}
\mathcal{L}(\hat{u},\hat{f})-\mathcal{L}(u^*,f^*)\le \underset{\bar{u}\in\mathcal{N}_u ,\bar{f}\in\mathcal{N}_f}{\inf}[\mathcal{L}(\bar{u},\bar{f})-\mathcal{L}(u^*,f^*)]+2\underset{u\in\mathcal{N}_u ,f\in\mathcal{N}_f}{\sup}\left |\hat{\mathcal{L}}(u,f)-\mathcal{L}(u,f)\right |.
\end{equation*}
\end{proof}

\subsection{Approximation Error}
\begin{assumption}\label{assumption:boundedness}(Boundedness). For simplicity, we would assume functions in $\mathcal{N}_{u}$ vanish on the boundary and\\
(i) $ \left\{f^{*}\right\} \cup \mathcal{N}_{f} \subseteq\left\{f:\left\|f ; L_{\infty}(\Omega)\right\| \leq B_{f}\right\}$ with $B_{f} \geq 1 $;\\
(ii)  $\left\{u^{*},u_d\right\} \cup \mathcal{N}_{u} \subseteq\left\{u\right.  :  \left.\left\|u ; W^{1,\infty}(\Omega)\right\| \leq B_{u}\right\}$ with $B_{u} \geq 1$;\\
(iii) $\left \|\Delta^{\alpha/2}u\right \|_ {L^2(\Omega)}\le B_0$ with $0<\alpha<1$.
\end{assumption} 

\begin{lemma}\label{lemma:laplace transform}(Proposition 3.6,\cite{dinezza2011hitchhikers})
Let $s\in(0,1)$ and $u\in H^{s}\left(\mathbb{R}^{n}\right)$, then 
\begin{equation*}
[u]_{H^{s}\left(\mathbb{R}^{n}\right)}^{2}=2 C(d, s)^{-1}\left\|(-\Delta)^{\frac{s}{2}} u\right\|_{L^{2}\left(\mathbb{R}^{n}\right)}^{2}
\end{equation*}
where C(d,s) is defined by (\ref{definition:C(d,s)}).
\end{lemma}

\begin{lemma} Based on the assumption \ref{assumption:boundedness}, for $\forall u\in N_u, f\in N_f$, we have:
\begin{equation*}
\mathcal{L}(u,f)-\mathcal{L}(u^*,f^*)\le 4B_u\left \| u-u^* \right \|_{L^{2}(\Omega)} + C(\Omega,\alpha){\left \| u-u^*\right \|^2_{H^{1}(\Omega)}} +(2B_0+2B_f)\left \| f-f^* \right \|_{L^{2}(\Omega)}.
\end{equation*}
\end{lemma}
\begin{proof} For $\forall u\in N_u, f\in N_f$,
\begin{equation*}
\mathcal{L}(u,f)-\mathcal{L}(u^*,f^*)=\mathcal{L}(u,f)-\mathcal{L}(u,f^*)+\mathcal{L}(u,f^*)-\mathcal{L}(u^*,f^*).
\end{equation*}
According to the definition of $\mathcal{L}$ in (\ref{eq:loss}),we have:
\begin{equation*}
\begin{aligned}
\mathcal{L}(u,f^*)-\mathcal{L}(u^*,f^*)&=(\| u-u_d\|_{L^2(\Omega)}^2 - \| u^*-u_d \|_{L^2(\Omega)}^2)+(\|\Delta^{\alpha/2}u+f^* \|_{L^2(\Omega)}^2\\
&-\| \Delta^{\alpha/2}u^*+f^*\|_{L^2(\Omega)}^2) +( \| u \|_{L^2(\partial \Omega)}^2-\| u^* \|_{L^2(\partial \Omega)}^2).
\end{aligned}
\end{equation*}
Since $u^*$ is the solution of (\ref{eq:problemsetup}) when the forcing term is $f^*$, then $\left \| \Delta^{\alpha/2}u^*+f^* \right \|_{L^2(\Omega)}^2=0$ and $\left \| u^* \right \|_{L^2(\partial \Omega)}^2=0$. Thus,
\begin{equation*}
\mathcal{L}(u,f^*)-\mathcal{L}(u^*,f^*)=(\left \| u-u_d \right \|_{L^2(\Omega)}^2 - \left \| u^*-u_d \right \|_{L^2(\Omega)}^2)+\left \| \Delta^{\alpha/2}u+f^* \right \|_{L^2(\Omega)}^2+\left \| u \right \|_{L^2(\partial \Omega)}^2.
\end{equation*}
Next from triangular inequality,
\begin{equation}\label{equ:app:1}
\begin{aligned}
\left \| u-u_d \right \|_{L^2(\Omega)}^2 - \left \| u^*-u_d \right \|_{L^2(\Omega)}^2&=(\left \| u-u_d \right \|_{L^2(\Omega)}+\left \| u^*-u_d \right \|_{L^2(\Omega)})(\left \| u-u_d \right \|_{L^2(\Omega)}-\left \| u^*-u_d \right \|_{L^2(\Omega)})\\
&\le(\left \| u-u_d \right \|_{L^2(\Omega)}+\left \| u^*-u_d \right \|_{L^2(\Omega)})\left \| u-u^* \right \|_{L^2(\Omega)},\\
\end{aligned}
\end{equation}
from the assumption \ref{assumption:boundedness}, we obtain $\exists B_u, B_f \ge 1$ s.t.
\begin{equation*}
\left \{ u^*,u_d \right \} \cup N_u \subseteq \left \{ u:\left \| u \right \|_{W^{1,\infty}} \le B_u \right \},
\left \{ f^*\right \} \cup N_f \subseteq \left \{ u:\left \| f \right \|_{L^{\infty }} \le B_f \right \},
\end{equation*}
thus
\begin{equation}\label{equ:app:2}
(\left \| u-u_d \right \|_{L^2(\Omega)}+\left \| u^*-u_d \right \|_{L^2(\Omega)})\le 4B_u.
\end{equation}
Since $\Omega \subset R^{d}$ is a bounded domain with lipschitz boundary and according to the assumption \ref{assumption:boundedness} and lemma \ref{lemma:laplace transform}, then we have:
\begin{equation}\label{equ:app:3}
\begin{aligned}
\left \| \Delta^{\alpha/2}u+f^* \right \|_{L^2(\Omega)}^2+\left \| u \right \|_{L^2(\partial \Omega)}^2 &= \left \| \Delta^{\alpha/2}(u-u^*) \right \|_{L^2(\Omega)}^2,\\
&= C(\Omega,\alpha)[u-u^*]^2_{H^\alpha(\Omega)},\\
&\le C(\Omega,\alpha){\left \| u-u^*\right \|^2_{H^{1}(\Omega)}}.
\end{aligned}
\end{equation}
Hence, from (\ref{equ:app:1})-(\ref{equ:app:3}),we can conclude:
\begin{equation*}
\mathcal{L}(u,f)-\mathcal{L}(u_f,f)\le 4B_u\left \| u-u^* \right \|_{L^{2}(\Omega)} + C(\Omega,\alpha){\left \| u-u^*\right \|^2_{H^{1}(\Omega)}}.
\end{equation*}
For the next part, from the assumption \ref{assumption:boundedness}, we have:
\begin{equation*}
\begin{aligned}
\mathcal{L}(u,f)-\mathcal{L}(u,f^*)
&=\left \| \Delta^{\alpha/2}u+f \right \|_{L^2(\Omega)}^2-\left \| \Delta^{\alpha/2}u+f^* \right \|_{L^2(\Omega)}^2,\\
&\le (2B_0+2B_f)\left \| f-f^* \right \|_{L^{2}(\Omega)}.
\end{aligned}
\end{equation*}
In conclusion,
\begin{equation*}
\mathcal{L}(u,f)-\mathcal{L}(u^*,f^*)\le 4B_u\left \| u-u^* \right \|_{L^{2}(\Omega)} + C(\Omega,\alpha){\left \| u-u^*\right \|^2_{H^{1}(\Omega)}} +(2B_0+2B_f)\left \| f-f^* \right \|_{L^{2}(\Omega)}.
\end{equation*}
\end{proof}

\begin{theorem}From~\cite{gühring2020approximation}, we have:
Suppose $\zeta>0$ is arbitrary, $k\in \left\{0,...,j\right\} $, $\Omega \subseteq  (0,1)^d, d\in\mathbb{N}$ and $n\ge k+1$, then there exist constants  L, C, depending on  d, n, k, p, $\zeta$  with the following properties: \\
for any $f^* \in \mathcal{F}_{n, d, p}$ , where
\begin{equation*}
\mathcal{F}_{n, d, p}:=\left\{f \in W^{n, p}\left((0,1)^{d}\right):\|f\|_{W^{n, p}\left((0,1)^{d}\right)} \leq 1\right\}
\end{equation*}
there is a Tanh neural network  $\Phi_{\varepsilon, f}$  with $d$-dimensional input and one-dimensional output, at most L layers and at most  
\begin{equation*}
\begin{aligned}
&(i) C(\varepsilon)^{-d/n} \qquad &&for \ k=0 \\
&(ii)C(\varepsilon )^{-d/(n-k-\zeta)} \qquad &&for \ k\ge1 
\end{aligned}
\end{equation*}
nonzero weights bounded in absolute value by 
\begin{equation*}
  C\epsilon^{-\theta} = C \varepsilon^{-2-\left(2\left(d / p+d+k+\zeta_{(k=2)}\right)+d / p+d\right) /\left(n-k-\zeta_{(k=2)}\right)}
\end{equation*}
such that its output $f$ satisfies:
$\left\|f-f^*\right\|_{W^{k, p}\left((0,1)^{d}\right)} \leq \varepsilon$
\end{theorem} 

\begin{corollary}\label{corrollary:app:f}
 Given any $f^* \in H^{1}(\Omega)$,  there exist a $\tanh$-network f with d-dimensional input and one-dimensional output, at most $Clog(d+2)$ layers and at most $C(d)\varepsilon^{-d}$ nonzero weights bounded in absolute value by $C(d)\varepsilon^{-2-\frac{9}{2}d}$, such that
\begin{equation*}
\left \| f-f^* \right \|_{L^{2}(\Omega)}\le \varepsilon
\end{equation*}
\end{corollary}
 
\begin{corollary}\label{corrollary:app:u}Given any $u^* \in H^{2}(\Omega)$, there exist a $\tanh$- network u with d-dimensional input and one-dimensional output, at most $Clog(d+2)$ layers and at most $C(d)\varepsilon^{-\frac{d}{1-\zeta}}$ nonzero weights bounded in absolute value by $C(d)\varepsilon^{-\frac{9d+8}{2-2\zeta}}$, such that
\begin{equation*}
\left \| u-u^* \right \|_{H^{1}(\Omega)}\le \varepsilon
\end{equation*}
\end{corollary}

\subsection{Statistical Error}
\begin{lemma}
Let $P_j =\left \{ r_{\epsilon j}, r_{oj}, \xi_j, {r_{\epsilon j}}', {r_{oj}}',{\xi_j}' \right \} $ and define $K(x,y) = x\cdot y$, we have
\begin{equation*}
\begin{aligned}
&\mathbb{E}_{\left\{X_{i}\right\}_{i=1}^{N},\left\{Y_{k}\right\}_{k=1}^{N_g},\left\{Z_{l}\right\}_{l=1}^{N_u},\left\{P_{j}\right\}_{j=1}^{m}} \sup _{u \in \mathcal{N}_u,f \in \mathcal{N}_f} \pm[\mathcal{L}(u,f)-\widehat{\mathcal{L}}(u,f)] \\
\leq &\sum_{k=1}^{3} \mathbb{E}_{\left\{X_{i}\right\}_{i=1}^{N},\left\{Y_{k}\right\}_{k=1}^{N_g},\left\{Z_{l}\right\}_{l=1}^{N_u},\left\{P_{j}\right\}_{j=1}^{m}} \sup _{u \in \mathcal{N}_u,f \in \mathcal{N}_f} \pm\left[\mathcal{L}_{k}(u,f)-\widehat{\mathcal{L}}_{k}(u,f)\right]
\end{aligned}
\end{equation*}
Where 
\begin{equation*}
\mathcal{L}_{1}(u, f)=|\Omega| \mathbb{E}_{X \sim U(\Omega)}(L(u(X))-f(X))^2
\end{equation*}
\begin{equation*}
\mathcal{L}_{2}(u)=|\partial \Omega| \mathbb{E}_{Y \sim U(\partial \Omega)} u(Y)^2
\end{equation*}
\begin{equation*}
\mathcal{L}_{3}(u)=|\Omega| \mathbb{E}_{Z \sim U(\Omega)}(u(Z)-u_d(Z))^2
\end{equation*}
and $\widehat{\mathcal{L}}_{2}(u)$, $\widehat{\mathcal{L}}_{3}(u)$ is the discrete version of $\mathcal{L}_{2}(u)$, $\mathcal{L}_{3}(u)$, for example,
\begin{equation*}
\widehat{\mathcal{L}}_{2}(u)=\frac{|\partial \Omega| }{N_u} \sum_{i=1}^{N_u}u^2
\end{equation*}
Meanwhile, 
\begin{equation*}
\widehat{\mathcal{L}}_{1}(u,f)=|\Omega|\frac{1}{N}\frac{1}{m}\sum_{i=1}^{N}\sum_{j=1}^{m} K(\mu(X_i,P_j),\eta(X_i,P_j))
\end{equation*}
where $K(\mu(X_i,P_j),\eta(X_i,P_j))=\mu(X_i,P_j)\eta(X_i,P_j)$
\end{lemma}
\begin{proof}
Direct result from triangle inequality.
\end{proof}

\begin{proposition}
W.L.O.G, We may let $m=N-1$ and $j\neq i$. Denote $G=(X,P)$, $W$ and $V$ be the function spaces of $\mu$ and $\eta$. 
Let $\mathbb{D} = \{X_1,...,X_N, P_1,...,P_N\}$ be the training data, then we obtain:
\begin{equation}\label{eq:loss1}
\begin{aligned}
&\mathbb{E}_{\left\{X_{i}\right\}_{i=1}^{N},\left\{P_{j}\right\}_{j=1}^{N}} \sup _{u \in \mathcal{N}_u,f \in \mathcal{N}_f} \pm[\mathcal{L}_{1}(u,f)-\widehat{\mathcal{L}}_{1}(u,f)]\\
&=|\Omega|\left.\mathbb{E}_{\mathbb{D}}\left[\sup _{\mu \in W, \eta \in V} \mid \sum_{i \neq j} \frac{1}{N(N-1)} K\left(\mu\left(X_{i}, P_{j}\right), \eta\left(X_{i}, P_{j}\right)\right)-E_{G}[K(\mu(G), \eta(G))]\mid\right.\right]
\end{aligned}
\end{equation}
\end{proposition}
\begin{proof}
From equation (\ref{equ:laplace decomposition}), (\ref{equ:decomposition:1}), (\ref{equ:decomposition:2}),we have 
\begin{equation*}
\begin{aligned}
&\mathbb{E}_{\left\{X_{i}\right\}_{i=1}^{N},\left\{P_{j}\right\}_{j=1}^{N}} \sup _{u \in \mathcal{N}_u,f \in \mathcal{N}_f} \pm[\mathcal{L}_{1}(u,f)-\widehat{\mathcal{L}}_{1}(u,f)]\\
&=\mathbb{E}_{\mathbb{D}} \sup _{u \in \mathcal{N}_u,f \in \mathcal{N}_f} \pm\left[|\Omega| \mathbb{E}_{X \sim U(\Omega)}(L(u)-f)^2-|\Omega|\frac{1}{N}\frac{1}{N-1}\sum_{i=1}^{N}\sum_{j=1,i\neq j}^{N} K(\mu(X_i,P_j),\eta(X_i,P_j))\right]\\
&=|\Omega|\left.\mathbb{E}_{\mathbb{D}}\left[\sup _{\mu \in W, \eta \in V} \mid \sum_{i \neq j} \frac{1}{N(N-1)} K\left(\mu\left(X_{i}, P_{j}\right), \eta\left(X_{i}, P_{j}\right)\right)-E_{G}[K(\mu(G), \eta(G))]\mid\right.\right]
\end{aligned}
\end{equation*}
\end{proof} 

\begin{lemma}
 For any $G=(x,P)$,  ${G}'=({x}',{P}')$, we define $\mathcal{F}=\mathcal{F}_{\mu, \eta} \cup\left(-\mathcal{F}_{\mu, \eta}\right)$, where
\begin{equation*}
\mathcal{F}_{\mu, \eta}=\left\{K_{\mu, \eta}\left(G, G^{\prime}\right)=\frac{K\left(\mu\left(x, P^{\prime}\right), \eta\left(x, P^{\prime}\right)\right)+K\left(\mu\left(x^{\prime}, P\right), \eta\left(x^{\prime}, P\right)\right)}{2}: \mu \in W, \eta \in V\right\} .
\end{equation*}
Then $K_{\mu, \eta}\left(G, G^{\prime}\right)=K_{\mu, \eta}\left(G^{\prime}, G\right)$. Next, we define $\bar{K}_{\mu, \eta}(G)=\mathbb{E}_{{G}'}\left[K_{\mu, \eta}(G, {G}')\right]$, for the right side of (\ref{eq:loss1}), we have:
\begin{equation*}
\begin{aligned}
&\left.\mathbb{E}_{\mathbb{D}}\left[\sup _{\mu \in W, \eta \in V} \mid \sum_{i \neq j} \frac{1}{N(N-1)} K\left(\mu\left(X_{i}, P_{j}\right), \eta\left(X_{i}, P_{j}\right)\right)-E_{G}[K(\mu(G),\eta(G))]\mid\right.\right] \\
&=\mathbb{E}_{\mathbb{D}}\left[\sup _{K_{\mu, \eta} \in \mathcal{F}} \sum_{i \neq j} \frac{1}{N(N-1)}\left(K_{\mu, \eta}\left(G_{i}, G_{j}\right)-\mathbb{E}\left[K_{\mu, \eta}\left(G_{i}, G_{j}\right)\right]\right)\right] \\
&\leq \mathbb{E}_{\mathbb{D}}\left[\sup _{K_{\mu, \eta} \in \mathcal{F}} \sum_{i \neq j} \frac{1}{N(N-1)}\left(K_{\mu, \eta}\left(G_{i}, G_{j}\right)-\bar{K}_{\mu, \eta}\left(G_{i}\right)-\bar{K}_{\mu, \eta}\left(G_{j}\right)+\mathbb{E}\left[K_{\mu, \eta}\left(G_{i}, G_{j}\right)\right]\right)\right] \\
&+2 \mathbb{E}_{\mathbb{D}}\left[\sup_{K_{\mu, \eta} \in \mathcal{F}} \sum_{i=1}^{N} \frac{1}{N}\left(\bar{K}_{\mu, \eta}\left(G_{i}\right)-\mathbb{E}\left[\bar{K}_{\mu, \eta}\left(G_{i}\right)\right]\right)\right]
\end{aligned}
\end{equation*}
\end{lemma}
\begin{definition}
Let $\left\{\sigma_{k}\right\}_{k=1}^{N}$  be  $N$  i.i.d Rademacher variables. Then the Rademacher complexity of function class  $\mathcal{N}$  associate with random sample  $\left\{X_{k}\right\}_{k=1}^{N}$  is defined as
\begin{equation*}
\mathfrak{R}_{N}(\mathcal{N})=\mathbb{E}_{\left\{X_{k}, \sigma_{k}\right\}_{k=1}^{N}}\left[\sup _{u \in \mathcal{N}} \frac{1}{N} \sum_{k=1}^{N} \sigma_{k} u\left(X_{k}\right)\right]    
\end{equation*}
\end{definition}

\begin{lemma}
\begin{equation*}
2 \mathbb{E}_{\mathbb{D}}\left[\sup _{K_{\mu, \eta} \in \mathcal{F}} \sum_{i=1}^{N} \frac{1}{N}\left(\bar{K}_{\mu, \eta}\left(G_{i}\right)-\mathbb{E}\left[\bar{K}_{\mu, \eta}\left(G_{i}\right)\right]\right)\right]\leq 4 \mathfrak R(\mathcal{F}_1)
\end{equation*}
where $\mathcal{F}_1=\left\{\bar{K}_{\mu, \eta}(G)=\mathbb{E}_{{G}'}\left[K_{\mu, \eta}(G, {G}')\right], K_{\mu, \eta}(G, {G}')\in \mathcal{F}\right\}$.
\end{lemma}
\begin{proof}
Take $\left\{\widetilde{G}_{i}\right\}_{i=1}^{N}$ as an independent copy of $\left\{G_{i}\right\}_{i=1}^{N}$, then
\begin{equation*}
\begin{aligned}
&2 \mathbb{E}_{\mathbb{D}}\left[\sup _{K_{\mu, \eta} \in \mathcal{F}} \sum_{i=1}^{N} \frac{1}{N}\left(\bar{K}_{\mu, \eta}\left(G_{i}\right)-\mathbb{E}\left[\bar{K}_{\mu, \eta}\left(G_{i}\right)\right)\right]\right]\\
&=\frac{2}{N} \mathbb{E}_{\mathbb{D}}\left[\sup _{K_{\mu, \eta} \in \mathcal{F}} \sum_{i=1}^{N} \left(\bar{K}_{\mu, \eta}\left(G_{i}\right)-\mathbb{E}\left[\bar{K}_{\mu, \eta}\left(\widetilde{G}_{i}\right)\right]\right)\right]\\
&\leq \frac{2}{N}\mathbb{E}_{\mathbb{D}} \mathbb{E}_{\mathbb{D'}}\left[\sup _{K_{\mu, \eta} \in \mathcal{F}} \sum_{i=1}^{N} \left(\bar{K}_{\mu, \eta}\left(G_{i}\right)-\bar{K}_{\mu, \eta}\left(\widetilde{G}_{i}\right)\right)\right]\\
&= \frac{2}{N}\mathbb{E}_{\mathbb{D}} \mathbb{E}_{\mathbb{D'}} \mathbb{E}_{\left\{\sigma_{i}\right\}_{i=1}^{N}}\left[\sup _{K_{\mu, \eta} \in \mathcal{F}} \sum_{i=1}^{N} \sigma_i\left(\bar{K}_{\mu, \eta}\left(G_{i}\right)-\bar{K}_{\mu, \eta}\left(\widetilde{G}_{i}\right)\right)\right]\\
&\leq \frac{2}{N}\mathbb{E}_{\mathbb{D}}\mathbb{E}_{\left\{\sigma_{i}\right\}_{i=1}^{N}}\left[\sup _{K_{\mu, \eta} \in \mathcal{F}} \sum_{i=1}^{N} \sigma_i\bar{K}_{\mu, \eta}\left(G_{i}\right)\right]+\frac{2}{N}\mathbb{E}_{\mathbb{D'}}\mathbb{E}_{\left\{\sigma_{i}\right\}_{i=1}^{N}}\left[\sup _{K_{\mu, \eta} \in \mathcal{F}} \sum_{i=1}^{N} -\sigma_i\bar{K}_{\mu, \eta}\left(\widetilde{G}_{i}\right)\right]\\
&=\frac{4}{N}\mathbb{E}_{\mathbb{D}}\mathbb{E}_{\left\{\sigma_{i}\right\}_{i=1}^{N}}\left[\sup _{K_{\mu, \eta} \in \mathcal{F}} \sum_{i=1}^{N} \sigma_i\bar{K}_{\mu, \eta}\left(G_{i}\right)\right]\\
&=4\mathfrak{R}_{N}(\mathcal{F}_1)
\end{aligned}
\end{equation*}
\end{proof}

For simplicity of notation, $\bar{K}_{\mu, \eta}(G)$ can also be denoted as $\bar{K}_{u, f}(G)$. Now we denote $\mathcal{F}_{u \times f}= {N}_u \times {N}_f$ and the metric in $\mathcal{F}_{u \times f}$ as: 
\begin{equation*}
 d_{\mathcal{F}_{u \times f}}((u,f),(\tilde{u},\tilde{f})) = \max \left \{  {\left \| u-\tilde{u} \right \|}_\infty, {\left \| f-\tilde{f} \right \|}_\infty \right \}.
\end{equation*}
By the boundedness of $u$ and $f$, we may have the following result:
\begin{lemma}\label{lemma:2:lipschitz}
If ${\left \| u \right \|}_\infty\leq B_u$ and ${\left \| f \right \|}_\infty\leq B_f$, then
\begin{equation*}
{\left \|\bar{K}_{u, f}(G)-\bar{K}_{\tilde{u}, \tilde{f}}(G)\right \|}_\infty \leq C(d,\alpha,r_0,\epsilon,B_u,B_f) \max \left \{  {\left \| u-\tilde{u} \right \|}_\infty, {\left \| f-\tilde{f} \right \|}_\infty \right \}   
\end{equation*}
\end{lemma}
\begin{proof}
Let $k_1= \frac{C_{d,\alpha }\left | S^{d-1} \right |{r_{0}}^{2-\alpha} }{2(2-\alpha)}$ and $k_2=\frac{C_{d,\alpha }\left | S^{d-1} \right |{r_0}^{-\alpha}}{2\alpha}$. If we completely expand the left side of above inequality, we can obtain:
\begin{equation*}
\begin{aligned}
&{\left \|\bar{K}_{u, f}(G)-\bar{K}_{\tilde{u}, \tilde{f}}(G)\right \|}_\infty = {\left \|\mathbb{E}_{{G_i}}\left[K_{u, f}(G, {G_i})\right]-\mathbb{E}_{{G_i}}\left[K_{\tilde{u}, \tilde{f}}(G, {G_i})\right]  \right \|}_\infty\\
\le & \| \frac{2k_1^2}{r_{\epsilon}^2 r_{\epsilon}'^2}(u(x)u(x)- \tilde{u}(x)\tilde{u}(x))- \frac{k_1 k_2}{r_{\epsilon}^2}(u(x)u(x-{r_{oi}}'{\xi_i}')-\tilde{u}(x)\tilde{u}(x-{r_{oi}}'{\xi_i}')) + ... \|_\infty 
\end{aligned}
\end{equation*}
Take $\frac{k_1 k_2}{r_{\epsilon}^2}(u(x)u(x-{r_{oi}}'{\xi_i}')-\tilde{u}(x)\tilde{u}(x-{r_{oi}}'{\xi_i}'))$ for example, we have 
\begin{equation*}
\begin{aligned}
&{\left \| \frac{k_1 k_2}{{r_{\epsilon}}^2 }(u(x)u(x-{r_{oi}}'{\xi_i}')-\tilde{u}(x)\tilde{u}(x-{r_{oi}}'{\xi_i}'))\right \|}_\infty\\
&={\left \|\frac{k_1 k_2}{{r_{\epsilon}}^2 }(u(x)u(x-{r_{oi}}'{\xi_i}')-u(x)\tilde{u}(x-{r_{oi}}'{\xi_i}')+u(x)\tilde{u}(x-{r_{oi}}'{\xi_i}')-\tilde{u}(x)\tilde{u}(x-{r_{oi}}'{\xi_i}'))\right \|}_\infty\\
&={\left \|\frac{k_1 k_2}{{r_{\epsilon}}^2 }(u(x)(u(x-{r_{oi}}'{\xi_i}')-\tilde{u}(x-{r_{oi}}'{\xi_i}'))+\tilde{u}(x-{r_{oi}}'{\xi_i}')(u(x)-\tilde{u}(x))\right \|}_\infty\\
&\leq \frac{k_1 k_2 B_u}{{\epsilon^2}} ({\left \| u(x)-\tilde{u}(x)  \right \|}_\infty+ {\left \| u(x-{r_{oi}}'{\xi_i}')-\tilde{u}(x-{r_{oi}}'{\xi_i}')\right \|}_\infty)
\end{aligned}
\end{equation*}
Then we obtain:
\begin{equation*}
\begin{aligned}
&{\left \|\bar{K}_{u, f}(G)-\bar{K}_{\tilde{u}, \tilde{f}}(G)\right \|}_\infty \\
&\leq L_1{\left \| u(x)-\tilde{u}(x)  \right \|}_\infty+L_2{\left \| u(x-{r_{oi}}'{\xi_i}')-\tilde{u}(x-{r_{oi}}'{\xi_i}')\right \|}_\infty+...+L_f {\left \| f(x_i)-\tilde{f}(x_i)  \right \|}_\infty\\
&\leq (L_1+L_2+...){\left \| u(x)-\tilde{u}(x)  \right \|}_\infty+L_f {\left \| f(x)-\tilde{f}(x)  \right \|}_\infty\\
&\leq C(d,\alpha,r_0,\epsilon,B_u,B_f) \max \left \{  {\left \| u-\tilde{u} \right \|}_\infty, {\left \| f-\tilde{f} \right \|}_\infty \right \}   
\end{aligned}
\end{equation*}
\end{proof}  
\begin{definition}\cite{Vershynin_2018}
An  $\epsilon$ -cover of a set  T  in a metric space  (S, $\tau$)  is a subset  $T_{c} \subset S $ such that for each  $t \in T$ , there exists a  $t_{c} \in T_{c}$ such that  $\tau\left(t, t_{c}\right) \leq \epsilon$ . The  $\epsilon$ -covering number of  T , denoted as  $\mathcal{C}(\epsilon, T, \tau$)  is defined to be the minimum cardinality among all  $\epsilon$ -cover of  T  with respect to the metric  $\tau$ .
\end{definition}
\begin{lemma}\label{lemma:dudley}(Dudley’s entropy formula~\cite{DUDLEY1967290})
Let  $\mathcal{F}$  be a function class and  $\|f\|_{\infty} \leq B$  for any  $f \in \mathcal{F}$ , we have
\begin{equation*}
\mathfrak{R}_{N}(\mathcal{F}) \leq \inf _{0<\delta<B / 2}\left(4 \delta+\frac{12}{\sqrt{N}} \int_{\delta}^{B / 2} \sqrt{\log \mathcal{C}\left(\epsilon, \mathcal{F},\|\cdot\|_{\infty}\right)} d \epsilon\right)
\end{equation*}
\end{lemma}
\begin{corollary}
For simplicity, we denote the lipschitz coefficient $C(d,\alpha,r_0,\epsilon,B_u,B_f)$ in Lemma \ref{lemma:2:lipschitz} as $C_1$. It is not difficult to find $\|\bar{K}_{u, f}(G)\|_{\infty} \leq B_{\bar{K}}= C(d,\alpha,r_0,\epsilon)(B_u+B_f)^2$ for any $\bar{K}_{u, f}(G) \in 
\mathcal{F}_1$ from the definition of $\bar{K}_{u,f}$. Since we have $ \left\{X_{i}\right\}_{i=1}^{N}$, then from Lemma \ref{lemma:dudley}, 
\begin{equation*}
\mathfrak{R}_{N}(\mathcal{F}_{1}) \leq \inf _{0<\delta<B_{\bar{K}}}\left(4 \delta+\frac{12}{\sqrt{N}} \int_{\delta}^{B_{\bar{K}}/2} \sqrt{\log \mathcal{C}\left(\frac{\epsilon}{C_1}, \mathcal{N}_u,\|\cdot\|_{\infty}\right)\mathcal{C}\left(\frac{\epsilon}{C_1}, \mathcal{N}_f,\|\cdot\|_{\infty}\right)} d \epsilon\right)
\end{equation*}
\end{corollary}
\begin{proof}
From Lemma \ref{lemma:2:lipschitz} we know the mapping $(u,f)\mapsto \bar{K}_{u, f} $ is $C_1$-Lipschitz, then we have
\begin{equation*}
\mathcal{C}\left(\epsilon, \mathcal{F}_1,\|\cdot\|_{\infty}\right)\le \mathcal{C}\left(\frac{\epsilon}{C_1}, \mathcal{F}_{u \times f},\|\cdot\|_{\infty}\right)
\end{equation*}
Accoding to the definition of covering number, it is not difficult to find 
\begin{equation*}
\mathcal{C}\left(\frac{\epsilon}{C_1}, \mathcal{F}_{u \times f},\|\cdot\|_{\infty}\right)\leq\mathcal{C}\left(\frac{\epsilon}{C_1},\mathcal{N}_u,\|\cdot\|_{\infty}\right)\mathcal{C}\left(\frac{\epsilon}{C_1}, \mathcal{N}_f,\|\cdot\|_{\infty}\right)
\end{equation*}
then we can obtain this corollary from Dudley's entropy formula.
\end{proof}
\begin{lemma}\cite{jiao2021error}\label{lemma:network:lipschitz}
Let  $\mathcal{D}, \mathfrak{n}_{\mathcal{D}}, n_{i} \in \mathbb{N}^{+}$, $n_{\mathcal{D}}=1$, $B_{\theta} \geq 1$  and  $\rho$  be a bounded Lipschitz continuous function with  $B_{\rho}, L_{\rho} \leq 1$ . Set the parameterized function class  $\mathcal{P}=\mathcal{N}_{\rho}\left(\mathcal{D}, \mathfrak{n}_{\mathcal{D}}, B_{\theta}\right)$ . For any  $f(x ; \theta) \in \mathcal{P}$, $f(x ; \theta)$  is  $\sqrt{\mathfrak{n}_{\mathcal{D}}} B_{\theta}^{\mathcal{D}-1}\left(\prod_{i=1}^{\mathcal{D}-1} n_{i}\right)$ -Lipschitz continuous with respect to variable  $\theta$ , i.e.,
\begin{equation*}
|f(x ; \theta)-f(x ; \widetilde{\theta})| \leq \sqrt{\mathfrak{n}_{\mathcal{D}}} B_{\theta}^{\mathcal{D}-1}\left(\prod_{i=1}^{\mathcal{D}-1} n_{i}\right)\|\theta-\widetilde{\theta}\|_{2}, \quad \forall x \in \Omega
\end{equation*}
\end{lemma} 
\begin{lemma}\label{lemma:network:lipschitz:covering}\cite{jiao2021error}Let  $\mathcal{F}$  be a parameterized class of functions:  $\mathcal{F}=\{f(x ; \theta): \theta \in \Theta\}$ . Let  $\|\cdot\|_{\Theta}$  be a norm on  $\Theta$  and let  $\|\cdot\|_{\mathcal{F}}$  be a norm on  $\mathcal{F}$ . Suppose that the mapping  $\theta \mapsto f(x ; \theta)$  is L-Lipschitz, that is,
\begin{equation*}
\|f(x ; \theta)-f(x ; \widetilde{\theta})\|_{\mathcal{F}} \leq L\|\theta-\widetilde{\theta}\|_{\Theta},
\end{equation*}
then for any  $\epsilon>0, \mathcal{C}\left(\epsilon, \mathcal{F},\|\cdot\|_{\mathcal{F}}\right) \leq \mathcal{C}\left(\epsilon / L, \Theta,\|\cdot\|_{\Theta}\right)$.
\end{lemma}

\begin{lemma}\label{lemma:network:covering}\cite{jiao2021error} Suppose that  $T \subset \mathbb{R}^{d}$  and  $\|t\|_{2} \leq B$  for  $t \in T$ , then
\begin{equation*}
\mathcal{C}\left(\epsilon, T,\|\cdot\|_{2}\right) \leq\left(\frac{2 B \sqrt{d}}{\epsilon}\right)^{d}
\end{equation*}
\end{lemma}
\begin{corollary}
Let$ B_{\theta} \geq 1, B_{\psi} \geq 1$, then 
\begin{equation*}
\mathcal{C}\left(\epsilon, \mathcal{N}_u,\|\cdot\|_{\infty}\right) \leq \left(\frac{2 B_{\theta}^{\mathcal{D}_u}\mathfrak{n}_{\mathcal{D}_u}\left(\prod_{i=1}^{\mathcal{D}_u-1} n_{i}\right)}{\epsilon} \right)^{\mathfrak{n}_{\mathcal{D}_u}}
\end{equation*}
\begin{equation*}
\mathcal{C}\left(\epsilon, \mathcal{N}_f,\|\cdot\|_{\infty}\right) \leq \left(\frac{2B_{\psi}^{\mathcal{D}_f}\mathfrak{n}_{\mathcal{D}_f}\left(\prod_{i=1}^{\mathcal{D}_f-1} n_{i}\right)}{\epsilon} \right)^{\mathfrak{n}_{\mathcal{D}_f}}
\end{equation*}
\end{corollary}
\begin{proof}
It is not difficult to deduce this result from Lemma \ref{lemma:network:lipschitz}, Lemma \ref{lemma:network:lipschitz:covering}, Lemma \ref{lemma:network:covering}.
\end{proof} 
\begin{definition}
(\cite{10.1162/NECO_a_00028},Empirical Rademacher Chaos Complexity)
Let  $\mathcal{F}$  be a class of functions on  $X \times X$  and  $\left\{\epsilon_{i}: i \in \mathbb{N}_{n}\right\}$
are independent Rademacher random variables. Furthermore, $\mathbf{x}=\left\{x_{i}: i \in \mathbb{N}_{n}\right\} $ are independent random variables distributed according to a distribution  $\mu$  on  X . The homogeneous Rademacher chaos process of order two, with respect to the Rademacher variable  $\varepsilon$ , is a random variable system defined by  
\begin{equation*}
  \left\{\hat{U}_{f}(\varepsilon)=\frac{1}{n} \sum_{i, j \in \mathbb{N}_{n}, i<j} \varepsilon_{i} \varepsilon_{j} f\left(x_{i}, x_{j}\right): f \in \mathcal{F}\right\}
\end{equation*}

We refer to the expectation of its suprema: $\hat{\mathcal{U}}_{n}(\mathcal{F})=\mathbb{E}_{\varepsilon}\left[\sup _{f \in \mathcal{F}}\left|\hat{U}_{f}(\varepsilon)\right|\right] $ as the empirical Rademacher chaos complexity over  $\mathcal{F}$.
\end{definition}
\begin{lemma}
Let $\mathbb{D}'=\{G_i'\}_{i=1}^N$ be the ghost samples and $\tau_i$ be the i.i.d  Radamacher random variables, then we have:
\begin{equation*}
\begin{aligned}
&\mathbb{E}_{\mathbb{D}}\left[\sup _{K_{\mu, \eta} \in \mathcal{F}} \sum_{i \neq j} \frac{1}{N(N-1)}\left(K_{\mu, \eta}\left(G_{i}, G_{j}\right)-\bar{K}_{\mu, \eta}\left(G_{i}\right)-\bar{K}_{\mu, \eta}\left(G_{j}\right)+\mathbb{E}\left[K_{\mu, \eta}\left(G_{i}, G_{j}\right)\right]\right)\right] \\
&\leq \frac{8}{N-1}\mathbb{E}_{\mathbb{D}}\left[\hat{\mathcal{U}}_{N}(\mathcal{F})\right]
\end{aligned}
\end{equation*}
\end{lemma}
\begin{proof}
\begin{equation*}
\begin{aligned}
&\mathbb{E}_{\mathbb{D}}\left[\sup _{K_{\mu, \eta} \in \mathcal{F}} \sum_{i \neq j} \frac{1}{N(N-1)}\left(K_{\mu, \eta}\left(G_{i}, G_{j}\right)-\bar{K}_{\mu, \eta}\left(G_{i}\right)-\bar{K}_{\mu, \eta}\left(G_{j}\right)+\mathbb{E}\left[K_{\mu, \eta}\left(G_{i}, G_{j}\right)\right]\right)\right] \\
&=\mathbb{E}_{\mathbb{D}}\left[\sup _{K_{\mu, \eta} \in \mathcal{F}} \sum_{i \neq j} \frac{1}{N(N-1)} K_{\mu, \eta}\left(G_{i}, G_{j}\right)-\mathbb{E}_{\mathbb{D}^{\prime}}\left[K_{\mu, \eta}\left(G_{i}, G_{j}^{\prime}\right)\right]-\mathbb{E}_{\mathbb{D}^{\prime}}\left[K_{\mu, \eta}\left(G_{i}^{\prime}, G_{j}\right)\right]+\mathbb{E}_{\mathbb{D}^{\prime}}\left[K_{\mu, \eta}\left(G_{i}^{\prime}, G_{j}^{\prime}\right)\right]\right] \\
&\leq \mathbb{E}_{\mathbb{D}, \mathbb{D}^{\prime}}\left[\sup _{K_{\mu, \eta} \in \mathcal{F}} \sum_{i \neq j} \frac{1}{N(N-1)}\left(K_{\mu, \eta}\left(G_{i}, G_{j}\right)-K_{\mu, \eta}\left(G_{i}, G_{j}^{\prime}\right)-K_{\mu, \eta}\left(G_{i}^{\prime}, G_{j}\right)+K_{\mu, \eta}\left(G_{i}^{\prime}, G_{j}^{\prime}\right)\right)\right] \\
&=\mathbb{E}_{\mathbb{D}, \mathbb{D}^{\prime}, \tau}\left[\sup _{K_{\mu, \eta} \in \mathcal{F}} \sum_{i \neq j} \frac{1}{N(N-1)} \tau_{i} \tau_{j}\left(K_{\mu, \eta}\left(G_{i}, G_{j}\right)-K_{\mu, \eta}\left(G_{i}, G_{j}^{\prime}\right)-K_{\mu, \eta}\left(G_{i}^{\prime}, G_{j}\right)+K_{\mu, \eta}\left(G_{i}^{\prime}, G_{j}^{\prime}\right)\right)\right]\\
&\leq 4 \mathbb{E}_{\mathbb{D},  \tau}\left[\sup_{K_{\mu, \eta} \in \mathcal{F}} |\sum_{i \neq  j} \frac{1}{N(N-1)} \tau_{i} \tau_{j}\left(K_{\mu, \eta}\left(G_{i}, G_{j}\right)\right)|\right]\\
&\leq \frac{8}{N-1}\mathbb{E}_{\mathbb{D},  \tau}\left[\sup_{K_{\mu, \eta} \in \mathcal{F}} |\sum_{i <  j} \frac{1}{N} \tau_{i} \tau_{j}\left(K_{\mu, \eta}\left(G_{i}, G_{j}\right)\right)|\right]\\
&=\frac{8}{N-1}\mathbb{E}_{\mathbb{D}}\left[\hat{\mathcal{U}}_{N}(\mathcal{F})\right]
\end{aligned}
\end{equation*}
\end{proof}
\begin{lemma}\label{lemma:chaos}\cite{10.1162/NECO_a_00028}
Let  $\left\{f_{\ell}: \ell \in \mathbb{N}_{N}\right\}$  be a finite class of functions on  $X \times X$  and  $\left\{\epsilon_{i}\right.  :  \left.i \in \mathbb{N}_{n}\right\}$  are independent Rademacher random variables. Consider the homogeneous Rademacher chaos process of order two  $\left\{\hat{U}_{f_{\ell}}(\varepsilon)=\frac{1}{n} \sum_{i, j \in \mathbb{N}_{n}, i<j} \varepsilon_{i} \varepsilon_{j} f_{\ell}\left(x_{i}, x_{j}\right): \ell \in\right.   \left.\mathbb{N}_{N}\right\}$ . Then, we have that
\begin{equation*}
\mathbb{E}\left[\max _{\ell \in \mathbb{N}_{N}}\left|\hat{U}_{f_{\ell}}(\varepsilon)\right|\right] \leq 2 e \log (1+N) \max _{\ell \in \mathbb{N}_{N}}\left(\frac{1}{n^{2}} \sum_{i<j}\left|f_{\ell}\left(x_{i}, x_{j}\right)\right|^{2}\right)^{\frac{1}{2}},
\end{equation*}
where  $\mathbb{E}[\cdot]$ denotes the expectation with respect to the Rademacher variable  $\varepsilon$ .
\end{lemma}
\begin{theorem}

From the definition of cover, for any $K_{\mu, \eta} \in \mathcal{F}$, we have $K_{{\mu}', {\eta}'}\in {\mathcal{F}}'$, where $\mathcal{F}'$ is a finite class of functions, such that $\|K_{\mu, \eta}-K_{{\mu}', {\eta}'}\|_\mathcal{F}\leq \delta$ and $|\mathcal{F}'|=\mathcal{C}\left(\delta, \mathcal{F},\|\cdot\|_{\mathcal{F}}\right)$. Considering the definition of $K_{\mu, \eta}$,  we can set  $ M=\max \left\{K_{{\mu'}, {\eta'}}, \forall K_{{\mu'}, {\eta}'}\in{\mathcal{F}}'\right\}$. Then for simplicity, let $C_1=C(d,\alpha,r_0,\epsilon,B_u,B_f)$, we have 
\begin{equation*}
\hat{\mathcal{U}}_{N}(\mathcal{F})\leq N\delta+ 2 e \log (1+\mathcal{C}\left(\frac{\delta}{C_1}, \mathcal{N}_u,\|\cdot\|_{\infty}\right)\mathcal{C}\left(\frac{\delta}{C_1}, \mathcal{N}_f,\|\cdot\|_{\infty}\right)) M
\end{equation*}
\end{theorem}
\begin{proof}
Similar as Lemma \ref{lemma:2:lipschitz}, it is easy to find that$K_{{\mu}, {\eta}}$is also Lipschitz continuous on (u,f), and the Lipschitz coefficient is still $C_1=C(d,\alpha,r_0,\epsilon,B_u,B_f)$. Next from Lemma \ref{lemma:chaos}, we can obtain 
\begin{equation*}
\begin{aligned}
\hat{\mathcal{U}}_{N}(\mathcal{F})&=\mathbb{E}_{\tau}\left[\sup _{K_{\mu,\eta} \in \mathcal{F}}\left|\hat{U}_{K_{\mu,\eta}}(\tau)\right|\right]\\
&=\mathbb{E}_{\tau}\left[\sup _{K_{\mu,\eta} \in \mathcal{F},K_{{\mu}', {\eta}'}\in {\mathcal{F}}'}\left|\hat{U}_{K_{\mu,\eta}}(\tau)-\hat{U}_{K_{{\mu}',{\eta}'}}(\tau)+\hat{U}_{K_{{\mu}',{\eta}'}}(\tau)\right|\right]\\
&\leq \mathbb{E}_{\tau}\left[\sup _{K_{\mu,\eta} \in \mathcal{F},K_{{\mu}', {\eta}'}\in {\mathcal{F}}'}\left|\hat{U}_{K_{\mu,\eta}}(\tau)-\hat{U}_{K_{{\mu}',{\eta}'}}(\tau)\right|\right]+\mathbb{E}_{\tau}\left[\sup _{K_{{\mu}', {\eta}'}\in {\mathcal{F}}'}\left|\hat{U}_{K_{{\mu}',{\eta}'}}(\tau)\right|\right]\\
&=\mathbb{E}_{\tau}\left[\sup _{K_{\mu,\eta} \in \mathcal{F},K_{{\mu}', {\eta}'}\in {\mathcal{F}}'}\left|\sum_{i <  j}\frac{1}{N} \tau_{i} \tau_{j}\left(K_{\mu, \eta}\left(G_{i}, G_{j}\right)-K_{{\mu}', {\eta}'}\left(G_{i}, G_{j}\right)\right)\right|\right]+\mathbb{E}_{\tau}\left[\sup _{K_{{\mu}', {\eta}'}\in {\mathcal{F}}'}\left|\hat{U}_{K_{{\mu}',{\eta}'}}(\tau)\right|\right]\\
&\leq N\delta+2e\log(1+\mathcal{C}\left(\delta, \mathcal{F},\|\cdot\|_{\mathcal{F}}\right))\max _{K_{{\mu}', {\eta}'}\in {\mathcal{F}}'}\left(\frac{1}{N^{2}} \sum_{i<j}\left|K_{{\mu}', {\eta}'}\left(G_{i}, G_{j}\right)\right|^{2}\right)^{\frac{1}{2}}\\
&\leq N\delta+ 2 e \log (1+\mathcal{C}\left(\frac{\delta}{C_1}, \mathcal{N}_u,\|\cdot\|_{\infty}\right)\mathcal{C}\left(\frac{\delta}{C_1}, \mathcal{N}_f,\|\cdot\|_{\infty}\right)) M
\end{aligned}
\end{equation*}
\end{proof}

\begin{theorem}\label{Thm:sta_err}
W.L.O.G, Let $N_g=N_u=N$. Let $\mathcal{D}=\max\left\{\mathcal{D}_u,\mathcal{D}_f\right\}$, $\mathfrak{n}=\max\left\{\mathfrak{n}_{\mathcal{D}_u},\mathfrak{n}_{\mathcal{D}_f}\right\}$, $\mathcal{B}=\max\left\{B_{\theta},B_{\psi}\right\}$, then we obtain:
\begin{equation*}
\begin{aligned}
&\mathbb{E}_{\left\{X_{i}\right\}_{i=1}^{N},\left\{Y_{k}\right\}_{j=1}^{N_g},\left\{Z_{l}\right\}_{l=1}^{N_u},\left\{P_{j}\right\}_{j=1}^{m}} \sup _{u \in \mathcal{N}_u,f \in \mathcal{N}_f} \pm[\mathcal{L}(u,f)-\widehat{\mathcal{L}}(u,f)] \\
&\leq C(\Omega,\alpha,d,r_0,\epsilon,B_u,B_f)\frac{\mathfrak{n}^{2} \sqrt{\mathfrak{n}} \mathcal{B}^{2}{\mathcal{D}}}{\sqrt{N}}\log \mathfrak{n} \mathcal{B}N
\end{aligned}
\end{equation*}
\end{theorem}

\begin{proof}
 For k=1, 
\begin{equation*}
\begin{aligned}
&\mathbb{E}_{\left\{X_{i}\right\}_{i=1}^{N},\left\{P_{j}\right\}_{j=1}^{N}} \sup _{u \in \mathcal{N}_u,f \in \mathcal{N}_f} \pm[\mathcal{L}_{1}(u,f)-\widehat{\mathcal{L}}_{1}(u,f)]\\
&\leq |\Omega| \mathbb{E}_{\mathbb{D}}\left[\sup _{K_{\mu, \eta} \in \mathcal{F}} \sum_{i \neq j} \frac{1}{N(N-1)}\left(K_{\mu, \eta}\left(G_{i}, G_{j}\right)-\bar{K}_{\mu, \eta}\left(G_{i}\right)-\bar{K}_{\mu, \eta}\left(G_{j}\right)+\mathbb{E}\left[K_{\mu, \eta}\left(G_{i}, G_{j}\right)\right]\right)\right] \\
&+|\Omega| 2 \mathbb{E}_{\mathbb{D}}\left[\sup _{K_{\mu, \eta} \in \mathcal{F}} \sum_{i=1}^{N} \frac{1}{N}\left(\bar{K}_{\mu, \eta}\left(G_{i}\right)-\mathbb{E}\left[\bar{K}_{\mu, \eta}\left(G_{i}\right)\right]\right)\right]
\end{aligned}
\end{equation*}

For simplicity, we denote $C(d,\alpha,r_0,\epsilon,B_u,B_f)$ as $C_1$, which refers to constant depending on $d,\alpha,r_0,\epsilon,B_u,B_f$. Now for part 2, we have
\begin{equation*}
\begin{aligned}
&2 \mathbb{E}_{\mathbb{D}}\left[\sup _{K_{\mu, \eta} \in \mathcal{F}} \sum_{i=1}^{N} \frac{1}{N}\left(\bar{K}_{\mu, \eta}\left(G_{i}\right)-\mathbb{E}\left[\bar{K}_{\mu, \eta}\left(G_{i}\right)\right]\right)\right]\leq 4 \mathfrak{R}(\mathcal{F}_{1})\\
& \leq 4\inf _{0<\delta<B_{\bar{K}}/2}\left(4 \delta+\frac{12}{\sqrt{N}} \int_{\delta}^{B_{\bar{K}}} \sqrt{\log \mathcal{C}\left(\frac{\epsilon}{C_1}, \mathcal{N}_u,\|\cdot\|_{\infty}\right)\mathcal{C}\left(\frac{\epsilon}{C_1}, \mathcal{N}_f,\|\cdot\|_{\infty}\right)} d \epsilon\right)\\
& \leq 4 \inf _{0<\delta< B_{\bar{K}}}\left(4 \delta+\frac{12}{\sqrt{N}} \int_{\delta}^{B_{\bar{K}}/2} \sqrt{\log \left(\frac{2 B_{\theta}^{\mathcal{D}_u}\mathfrak{n}_{\mathcal{D}_u}\left(\prod_{i=1}^{\mathcal{D}_u-1} n_{i}\right)}{\frac{\epsilon}{C_1}} \right)^{\mathfrak{n}_{\mathcal{D}_u}}\left(\frac{2B_{\psi}^{\mathcal{D}_f}\mathfrak{n}_{\mathcal{D}_f}\left(\prod_{i=1}^{\mathcal{D}_f-1} n_{i}\right)}{\frac{\epsilon}{C_1}} \right)^{\mathfrak{n}_{\mathcal{D}_f}}} d \epsilon\right)
\end{aligned}
\end{equation*}
Choosing $\delta=\frac{1}{\sqrt{N}}$, we next obtain
\begin{equation*}
\begin{aligned}
&2 \mathbb{E}_{\mathbb{D}}\left[\sup _{K_{\mu, \eta} \in \mathcal{F}} \sum_{i=1}^{N} \frac{1}{N}\left(\bar{K}_{\mu, \eta}\left(G_{i}\right)-\mathbb{E}\left[\bar{K}_{\mu, \eta}\left(G_{i}\right)\right]\right)\right]\\
&\leq 4C_1 \frac{(\sqrt{\mathfrak{n}_{\mathcal{D}_u}}+ \sqrt{\mathfrak{n}_{\mathcal{D}_f}})(\sqrt{\mathcal{D}_u}+\sqrt{\mathcal{D}_f})}{\sqrt{N}}\sqrt{\log \mathfrak{n}_{\mathcal{D}_u}\mathfrak{n}_{\mathcal{D}_f} \mathcal{B}_{\theta}\mathcal{B}_{\psi}N} \\
&\leq 4C_1 \frac{\sqrt{\mathfrak{n}} \sqrt{\mathcal{D}}}{\sqrt{N}}\sqrt{\log(\mathfrak{n} \mathcal{B}N)} 
\end{aligned}
\end{equation*}
Similarly for part 1, 
\begin{equation*}
\begin{aligned}
&\mathbb{E}_{\mathbb{D}}\left[\sup _{K_{\mu, \eta} \in \mathcal{F}} \sum_{i \neq j} \frac{1}{N(N-1)}\left(K_{\mu, \eta}\left(G_{i}, G_{j}\right)-\bar{K}_{\mu, \eta}\left(G_{i}\right)-\bar{K}_{\mu, \eta}\left(G_{j}\right)+\mathbb{E}\left[K_{\mu, \eta}\left(G_{i}, G_{j}\right)\right]\right)\right] \\
&\leq \frac{8}{N-1}(N\delta+ 2 e \log (1+\mathcal{C}\left(\frac{\delta}{C_1},\mathcal{N}_u,\|\cdot\|_{\infty}\right)\mathcal{C}\left(\frac{\delta}{C_1}, \mathcal{N}_f,\|\cdot\|_{\infty}\right)) M)\\
& \leq \frac{8}{N-1}(N\delta+ 2 e \log (e\left(\frac{2 B_{\theta}^{\mathcal{D}_u}\mathfrak{n}_{\mathcal{D}_u}\left(\prod_{i=1}^{\mathcal{D}_u-1} n_{i}\right)}{\frac{\delta}{C_1}} \right)^{\mathfrak{n}_{\mathcal{D}_u}}\left(\frac{2 B_{\psi}^{\mathcal{D}_f}\mathfrak{n}_{\mathcal{D}_f}\left(\prod_{i=1}^{\mathcal{D}_f-1} n_{i}\right)}{\frac{\delta}{C_1}} \right)^{\mathfrak{n}_{\mathcal{D}_f}})M)
\end{aligned}
\end{equation*}
According to the definition of $K_{\mu, \eta}$, we know that M can be expressed by $ d,\alpha, r_0, \epsilon, B_u, B_f$. Choosing $\delta=\frac{1}{N}$, we next obtain
\begin{equation*}
\begin{aligned}
&\mathbb{E}_{\mathbb{D}}\left[\sup _{K_{\mu, \eta} \in \mathcal{F}} \sum_{i \neq j} \frac{1}{N(N-1)}\left(K_{\mu, \eta}\left(G_{i}, G_{j}\right)-\bar{K}_{\mu, \eta}\left(G_{i}\right)-\bar{K}_{\mu, \eta}\left(G_{j}\right)+\mathbb{E}\left[K_{\mu, \eta}\left(G_{i}, G_{j}\right)\right]\right)\right] \\
&\leq C_1(\frac{1}{N-1}+\frac{(\mathfrak{n}_{\mathcal{D}_u}+\mathfrak{n}_{\mathcal{D}})(\mathcal{D}_u+\mathcal{D}_f)}{N-1}\log (\mathfrak{n}_{\mathcal{D}_u}\mathfrak{n}_{\mathcal{D}_f} \mathcal{B}_{\theta}\mathcal{B}_{\psi}N)\\
&\leq C_1\frac{\mathfrak{n}\mathcal{D}}{N-1}\log \mathfrak{n} \mathcal{B}N
\end{aligned}
\end{equation*}
For k=2 and k=3, it is not difficult to obtain conclusions from \cite{jiao2021error},
\begin{equation*}
\begin{aligned}
&\mathbb{E}_{\left\{Y_{k}\right\}_{k=1}^{N_g}} \sup _{u \in \mathcal{N}_u}[\mathcal{L}_{2}(u)-\widehat{\mathcal{L}}_{2}(u)]\\
&\leq C(\partial\Omega)\frac{{\mathfrak{n}^{2}_{\mathcal{D}_u}} \sqrt{\mathfrak{n}_{\mathcal{D}_u}} {\mathcal{B}_\theta}^{2}\sqrt{\mathcal{D}_u}}{\sqrt{N}}\sqrt{\log \mathfrak{n}_{\mathcal{D}_u} \mathcal{B_\theta}N} \\
&\leq C(\partial\Omega)\frac{\mathfrak{n}^{2} \sqrt{\mathfrak{n}} \mathcal{B}^{2}\sqrt{\mathcal{D}}}{\sqrt{N}}\sqrt{\log \mathfrak{n} \mathcal{B}N} 
\end{aligned}
\end{equation*}
Then we can also find:
\begin{equation*}
\begin{aligned}
&\mathbb{E}_{\left\{Z_{l}\right\}_{l=1}^{N_u}} \sup _{u \in \mathcal{N}_u}[\mathcal{L}_{3}(u)-\widehat{\mathcal{L}}_{3}(u)]\\
&\leq C(\Omega,B_u)\frac{\mathfrak{n}^{2} \sqrt{\mathfrak{n}} \mathcal{B}^{2}\sqrt{\mathcal{D}}}{\sqrt{N}}\sqrt{\log \mathfrak{n} \mathcal{B}N} 
\end{aligned}
\end{equation*}
Combining the results above, it is not difficult to obtain the conclusion.
\end{proof}

\subsection{Covergence Rate for the Inverse MC-fPINNs method}\label{covergence rate}
\begin{theorem}\label{theorem:convergence-rate}
Let $f^* \in H^{1}(\Omega)$ be the forcing term function and $u^* \in H^{2}(\Omega)$ be the corresponding solution of problem \eqref{eq:problemsetup}. Suppose the assumption \ref{assumption:boundedness} is fulfilled. Let $\rho$ be the $\tanh$ function $\frac{e^{x}-e^{-x}}{e^{x}+e^{-x}}$. Given any $\epsilon > 0$ and arbitrary $0<\zeta<1$, there exists a neural network $\hat{u}\in \mathcal{N}_{\rho}\left(\mathcal{D}_u, \mathfrak{n}_{\mathcal{D}_u}, B_{\theta}\right)$ and a neural network $\hat{f} \in \mathcal{N}_{\rho}\left(\mathcal{D}_f, \mathfrak{n}_{\mathcal{D}_f}, B_{\psi}\right)$, which satisfy 
\begin{equation*}
    \mathcal{D}_u,\mathcal{D}_f=Clog(d+1)
\end{equation*}
\begin{equation*}
\mathfrak{n}_{\mathcal{D}_u},\mathfrak{n}_{\mathcal{D}_f}= C(d)\varepsilon^{-\frac{d}{1-\zeta}}
\end{equation*}
\begin{equation*}
B_{\theta},B_{\psi}=C(d)\varepsilon^{-\frac{9d+8}{2-2\zeta}}
\end{equation*}
and number of samples:
\begin{equation*}
N=N_g=N_u=C(\Omega,\alpha,d,r_0,\epsilon,B_u,B_f) \varepsilon^{-\frac{23d+18-2\zeta}{1-\zeta}}
\end{equation*}
such that for $\forall\bar{u}\in\mathcal{N}_u,\forall\bar{f}\in\mathcal{N}_f$,we have:
\begin{equation*}
\mathcal{L}(\hat{u},\hat{f})-\mathcal{L}(u^*,f^*) < \varepsilon
\end{equation*}
\end{theorem}
\begin{proof}
According to Corollary \ref{corrollary:app:f} and \ref{corrollary:app:u}, there exit $\mathfrak{n}_{\mathcal{D}_u}=C(d)\varepsilon^{-\frac{d}{1-\zeta}}>C(d)\varepsilon^{-d}=\mathfrak{n}_{\mathcal{D}_f}$, $B_{\theta}=C(d)\varepsilon^{-\frac{9d+8}{2-2\mu}}>C(d)\varepsilon^{-2-\frac{9}{2}d}=B_{\psi}$ such that
\begin{equation*}
\left \| u-u^* \right \|_{H^{1}(\Omega)} < \varepsilon
\end{equation*}
\begin{equation*}
\left \| f-f^* \right \|_{L^{2}(\Omega)} <\varepsilon
\end{equation*}
In this case, we have 
\begin{equation}\label{equ:cov:app}
\mathcal{E}_{a p p}=\underset{\bar{u}\in\mathcal{N}_u ,\bar{f}\in\mathcal{N}_f}{\inf}[\mathcal{L}(\bar{u},\bar{f})-\mathcal{L}(u^*,f^*)]<C(\Omega,\alpha,B_u,B_f,B_0)\varepsilon   
\end{equation}
By appling Theorem \ref{Thm:sta_err} with $\mathcal{B}=C(d)\varepsilon^{-\frac{9d+8}{2-2\zeta}}, \mathfrak{n}=C(d)\varepsilon^{-\frac{d}{1-\zeta}}, D=Clog(d+2)$. Now we can conclude by setting 
\begin{equation*}
N=C(\Omega,\alpha,d,r_0,\epsilon,B_u,B_f) \varepsilon^{-(\frac{23d+18-2\zeta}{1-\zeta}+\tau)}
\end{equation*}
for an arbitrary small number $\tau>0$, we have 
\begin{equation}\label{equ:cov:sta}
\mathcal{E}_{sta}=\mathbb{E}_{\left\{X_{i}\right\}_{i=1}^{N},\left\{Y_{k}\right\}_{j=1}^{N_g},\left\{Z_{l}\right\}_{l=1}^{N_u},\left\{P_{j}\right\}_{j=1}^{m}} \sup _{u \in \mathcal{N}_u,f \in \mathcal{N}_f} \pm[\mathcal{L}(u,f)-\widehat{\mathcal{L}}(u,f)]<\varepsilon 
\end{equation}
Combining Proposition \ref{proposition:error-decomposition}, \eqref{equ:cov:app} and \eqref{equ:cov:sta} yields the conclusion.
\end{proof}

\section{Experiment}\label{experiment}
In this section, we will use the MC-fPINNs method to solve the inverse problem for the following Poisson equations:
\begin{equation*}
\begin{aligned}
(-\Delta)^{\alpha/2}u(x) ={f}(x),  \qquad  &x\in \Omega =\left\{x\mid {\left \| x  \right \|}_{2}^{2} \le  1 \right\},\\
u(x)=0,  \qquad  &x\in R^d \backslash\Omega,
\end{aligned}
\end{equation*}
with the exact solution being set as $u^{*}(x)=(1-{\left \| x  \right \|}^{2})_{+}^{1+\frac{\alpha}{2}}$, where $x_{+}:=max\{x,0\}$. By substituting $u^*$ in to the first equation, we can obtain the expression of the forcing term:
\begin{equation*}
f^{*}(x)= 2^\alpha \Gamma(\frac{\alpha}{2} +2)\Gamma (\frac{\alpha+d}{2} )\Gamma(\frac{d}{2})^{-1}(1-(1+\frac{\alpha}{d})\left \| x \right \|_2^2)
\end{equation*}
which will be used to evaluate the numerical error. The extra measurements of $u$ are set by applying a perturbation with a relative noise level of $\delta$:
\begin{equation*}
u_d^{\delta}(x)=u^{*}(x)+\delta u^{*}(x)\xi(x),  \qquad  x\in \Omega,
\end{equation*}
where $\xi(x)$ is a random variable following a standard Gaussian distribution.
In the following computations, we choose 1000 test points uniformly in the physical domain to compute the relative error for $\hat{u}$ and $\hat{f}$, which are defined as:
\begin{equation*}
\begin{aligned}
R e_{u} & :=\frac{\left\|\hat{u}(x)-u^*(x)\right\|_{L^{2}\left(\Omega\right)}}{\left\|u^*(x)\right\|_{L^{2}\left(\Omega\right)}} \\
R e_{f} & :=\frac{\left\|\hat{f}(x)-f^*(x)\right\|_{L^{2}(\Omega)}}{\left\|f^*(x)\right\|_{L^{2}(\Omega)}}
\end{aligned}
\end{equation*}
We set $\epsilon=0.01$, $r_0=0.3$, $m=30$ and construct two neural networks, each consisting of 4 hidden layers with 64 neurons per hidden layer, to approximate $u(x)$ and $f(x)$ respectively. The number of training epochs is set to be $10^4$ and the training points and test points are uniformly sampled in the physical domain. During the iteration, we use 256 residual points for each mini-batch training. As to the optimization, we use the Adam optimizer with a changing learning rate starting from $10^{-3}$ and $10^{-4}$ respectively for $ u_{NN}(x;\theta )$ and $f_{NN}(x;\psi)$. Furthermore, we would set the value of $u$ outside $\Omega$ as $0$ (i.e., we define $u(x) = \left(1- \|x\|_{2}^{2}\right)_{+} \cdot u_{N N}(x;\theta)$), thus there is no need to sample outside the domain $\Omega$.

\begin{table}
\centering
\begin{tabular}{cccc}
\toprule
&\multicolumn{3}{c}{$\delta$}   \\
\cmidrule(lr){2-4}
$\alpha$ & $1\%$ & $5\%$ & $10\%$  \\
\cmidrule{1-4}
0.5 & $8.54 \times 10^{-3}$ & $9.51 \times 10^{-3}$ &  $1.13 \times 10^{-2}$ \\
1.2 & $2.54 \times 10^{-2}$ & $3.94 \times 10^{-2}$ & $3.11 \times 10^{-2}$   \\
1.5 & $2.25 \times 10^{-2}$ & $2.91 \times 10^{-2}$ & $3.17 \times 10^{-2}$   \\
1.8 & $4.06 \times 10^{-2}$ & $3.70 \times 10^{-2}$ & $5.98 \times 10^{-2}$   \\
\bottomrule
\end{tabular}
\caption{The relative $L^{2}$-error of the reconstructed $\hat{f}$ with different $\alpha$ at different noise level}
\label{Table 1}
\end{table}

\begin{figure}
\centering
\subfloat[$\hat{f}$ with $\delta=1\%$.]
{\includegraphics[width=0.30\linewidth]{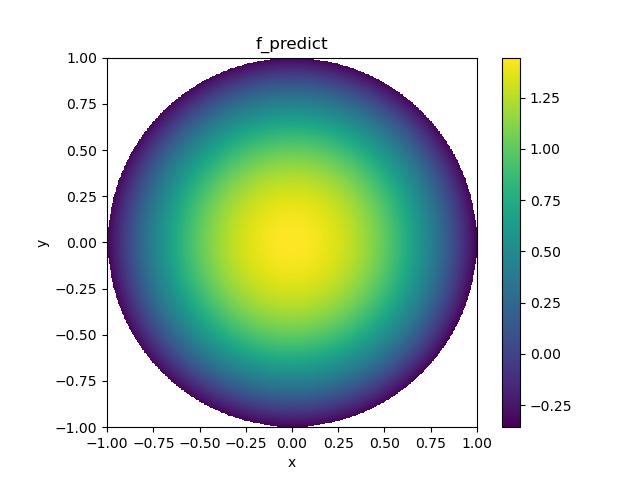}}
\subfloat[$\hat{f}$ with $\delta=5\%$.]
{\includegraphics[width=0.30\linewidth]{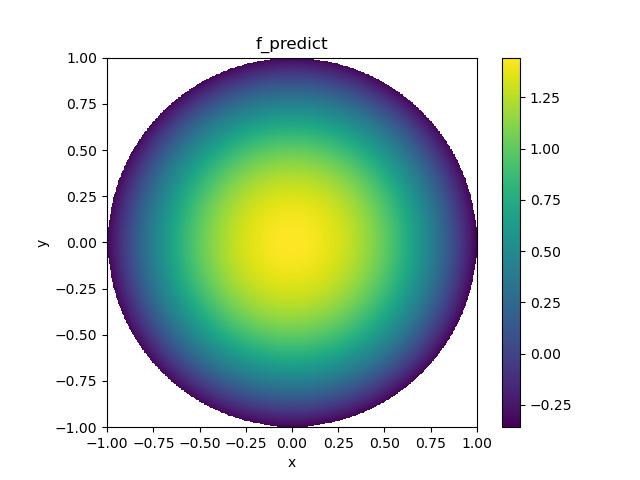}}
\subfloat[$\hat{f}$ with $\delta=10\%$.]
{\includegraphics[width=0.30\linewidth]{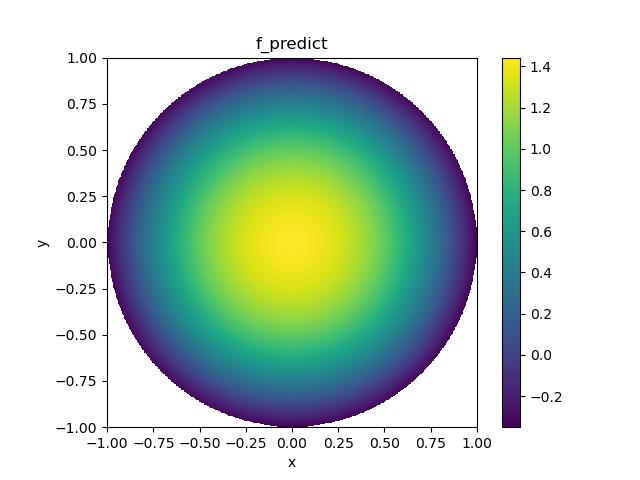}}
\\
\subfloat[$|\hat{f}-f^{*}|$ with $\delta=1\%$.]
{\includegraphics[width=0.30\linewidth]{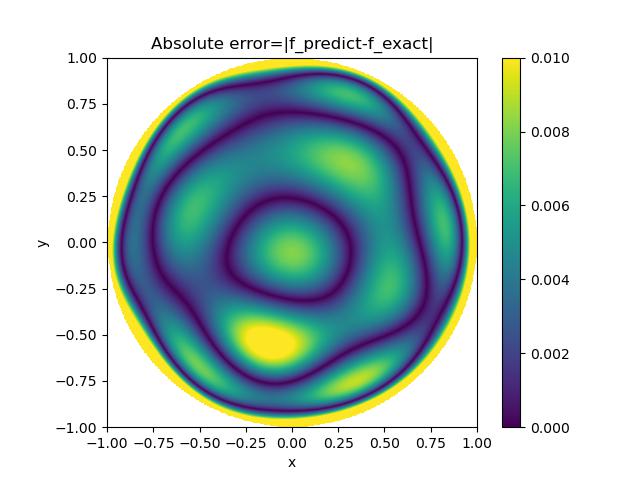}}
\subfloat[$|\hat{f}-f^{*}|$ with $\delta=5\%$.]
{\includegraphics[width=0.30\linewidth]{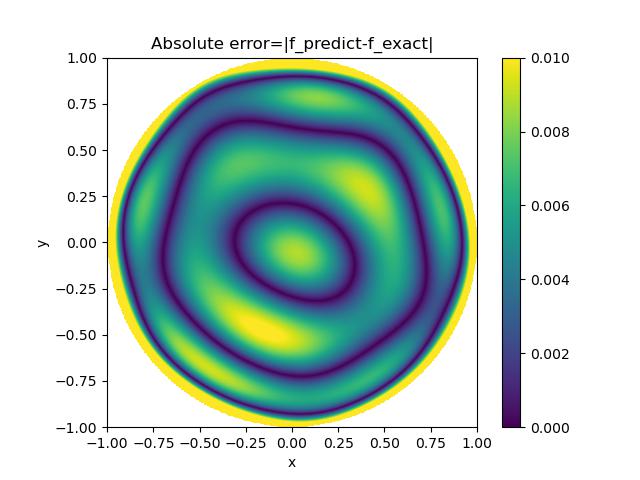}}
\subfloat[$|\hat{f}-f^{*}|$ with $\delta=10\%$.]
{\includegraphics[width=0.30\linewidth]{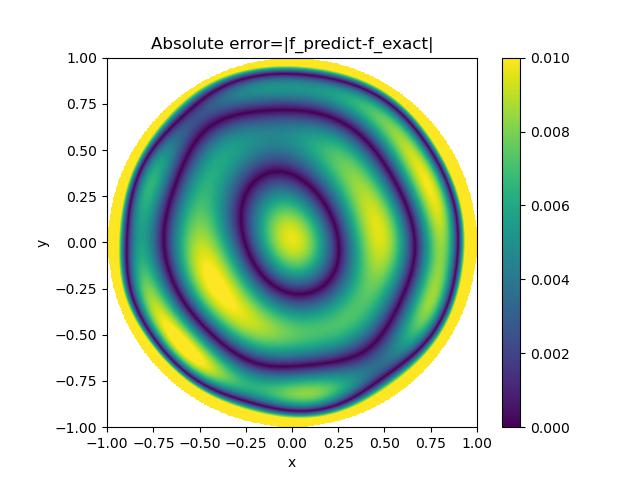}}
\caption{The reconstructions of $f$ (top) and the corresponding point-wise absolute error $|\hat{f}-f^{*}|$ (bottom) when $\alpha=0.5$}
\label{Figure 1}
\end{figure}

\begin{figure}
\centering
\subfloat[$\hat{f}$ with $\delta=1\%$.]
{\includegraphics[width=0.30\linewidth]{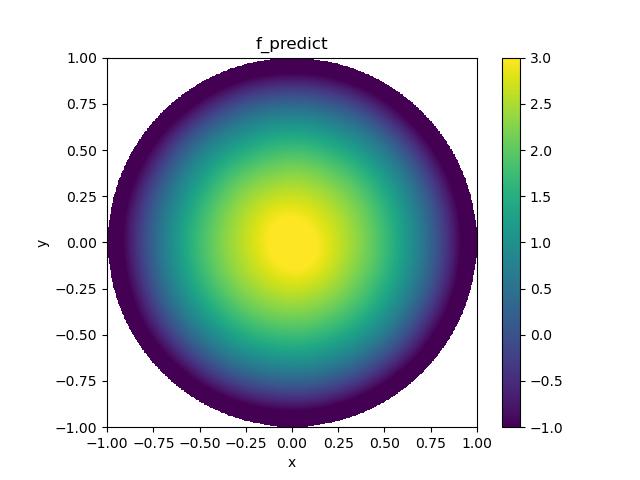}}
\subfloat[$\hat{f}$ with $\delta=5\%$.]
{\includegraphics[width=0.30\linewidth]{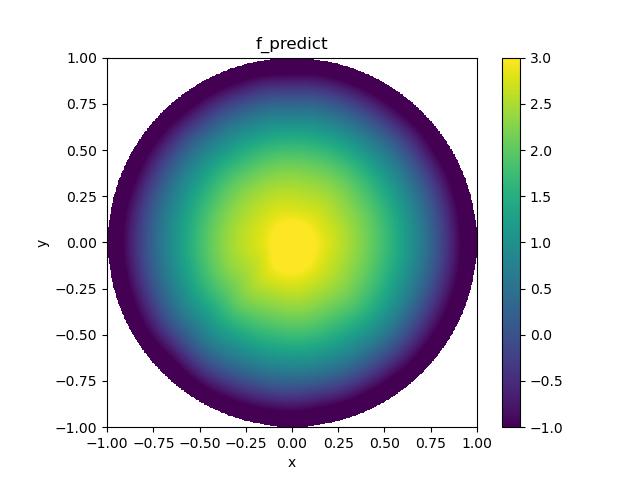}}
\subfloat[$\hat{f}$ with $\delta=10\%$.]
{\includegraphics[width=0.30\linewidth]{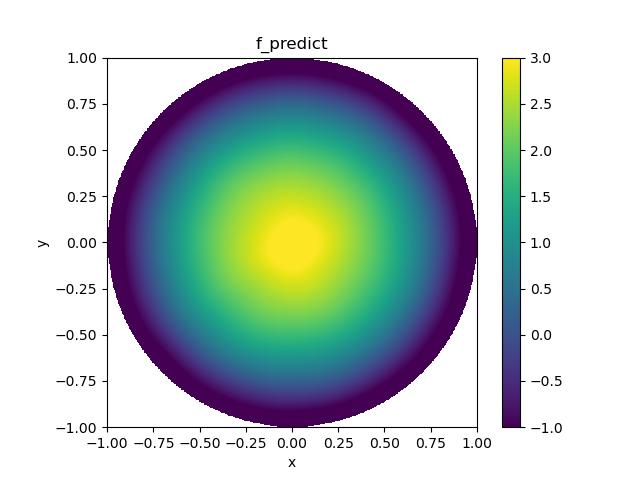}}
\\
\subfloat[$|\hat{f}-f^{*}|$ with $\delta=1\%$.]
{\includegraphics[width=0.30\linewidth]{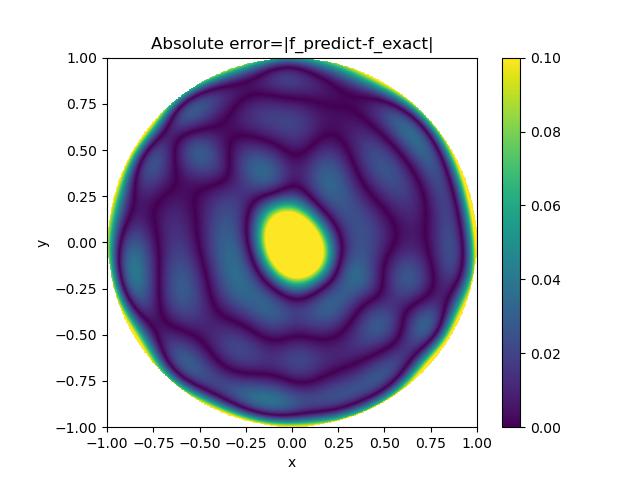}}
\subfloat[$|\hat{f}-f^{*}|$ with $\delta=5\%$.]
{\includegraphics[width=0.30\linewidth]{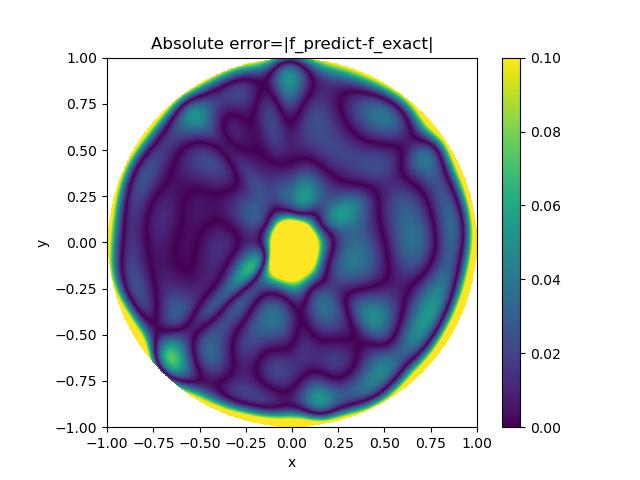}}
\subfloat[$|\hat{f}-f^{*}|$ with $\delta=10\%$.]
{\includegraphics[width=0.30\linewidth]{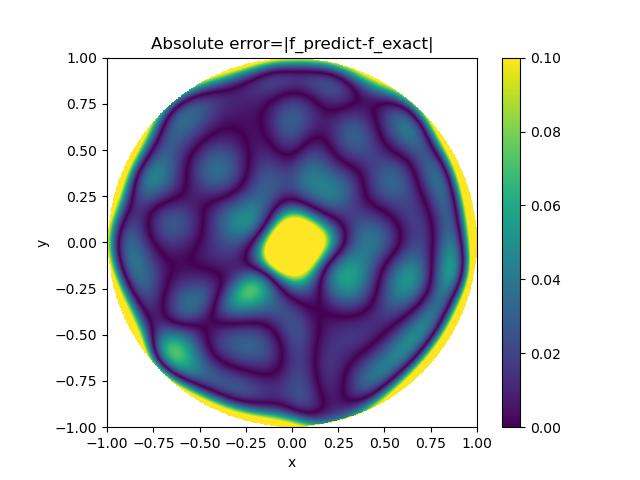}}
\caption{The reconstructions of $f$ (top) and the corresponding point-wise absolute error $|\hat{f}-f^{*}|$ (bottom) when $\alpha=1.2$}
\label{Figure 2}
\end{figure}

\begin{figure}
\centering
\subfloat[$\hat{f}$ with $\delta=1\%$.]
{\includegraphics[width=0.30\linewidth]{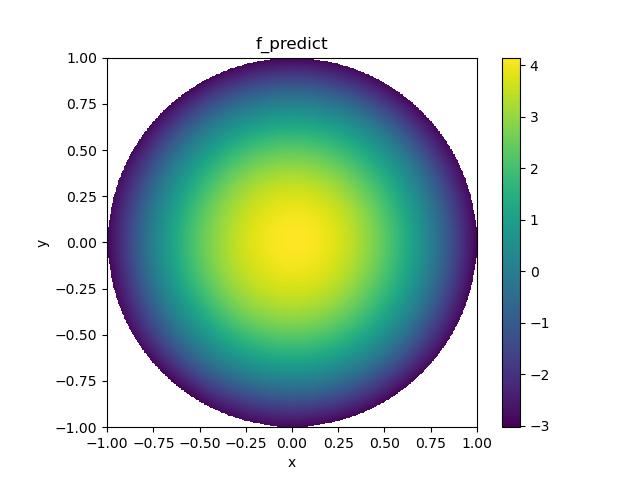}}
\subfloat[$\hat{f}$ with $\delta=5\%$.]
{\includegraphics[width=0.30\linewidth]{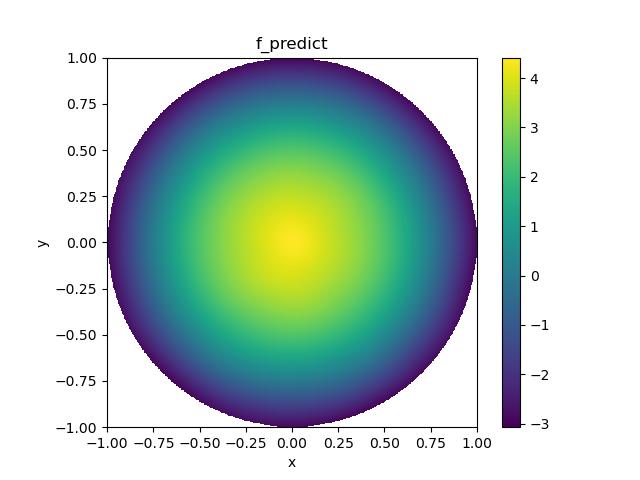}}
\subfloat[$\hat{f}$ with $\delta=10\%$.]
{\includegraphics[width=0.30\linewidth]{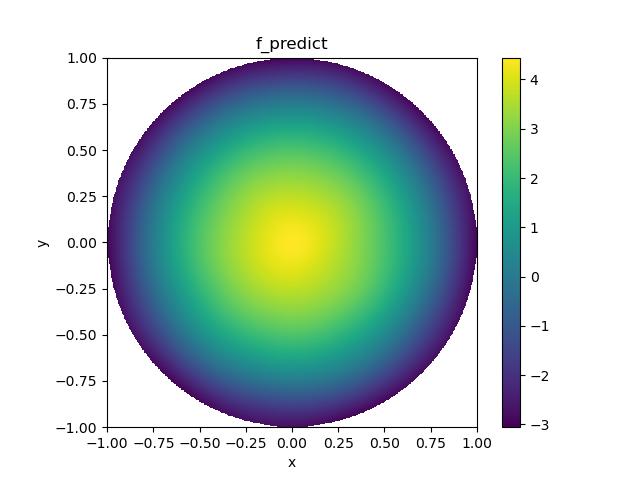}}
\\
\subfloat[$|\hat{f}-f^{*}|$ with $\delta=1\%$.]
{\includegraphics[width=0.30\linewidth]{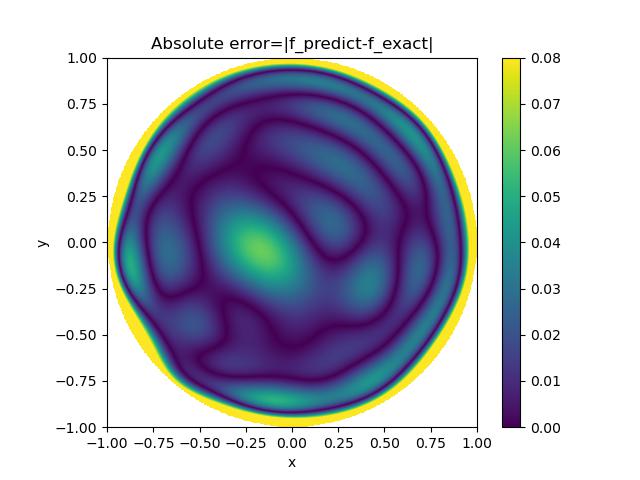}}
\subfloat[$|\hat{f}-f^{*}|$ with $\delta=5\%$.]
{\includegraphics[width=0.30\linewidth]{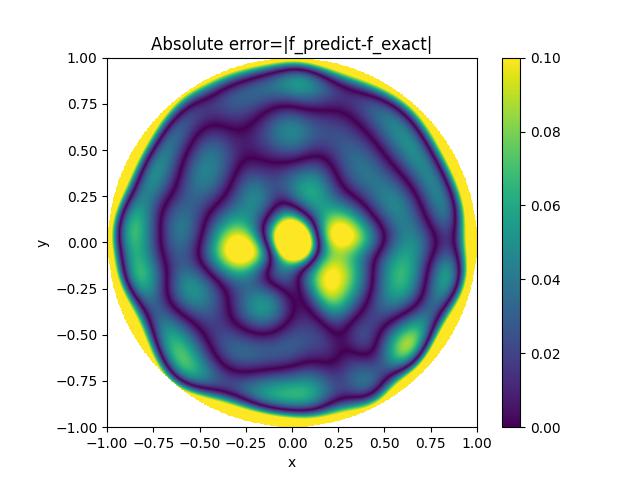}}
\subfloat[$|\hat{f}-f^{*}|$ with $\delta=10\%$.]
{\includegraphics[width=0.30\linewidth]{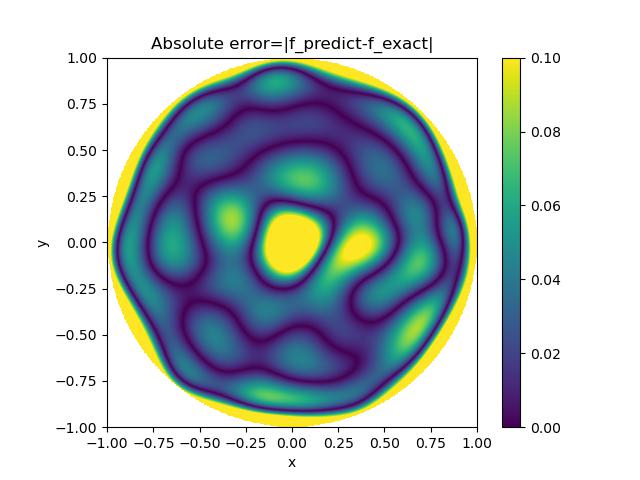}}
\caption{The reconstructions of $f$ (top) and the corresponding point-wise absolute error $|\hat{f}-f^{*}|$ (bottom) when $\alpha=1.5$}
\label{Figure 3}
\end{figure}

\begin{figure}
\centering
\subfloat[$\hat{f}$ with $\delta=1\%$.]
{\includegraphics[width=0.30\linewidth]{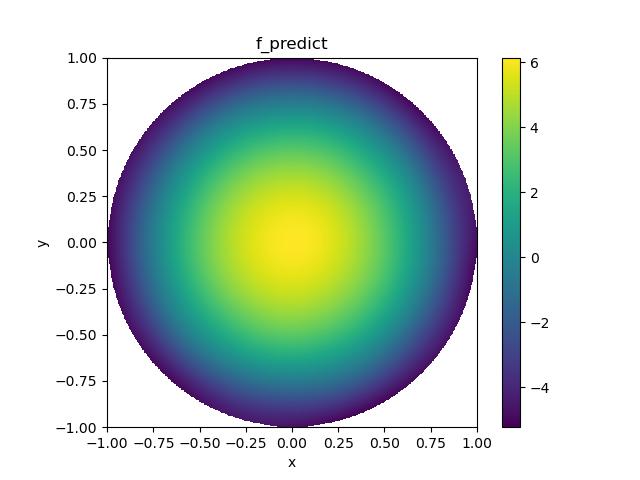}}
\subfloat[$\hat{f}$ with $\delta=5\%$.]
{\includegraphics[width=0.30\linewidth]{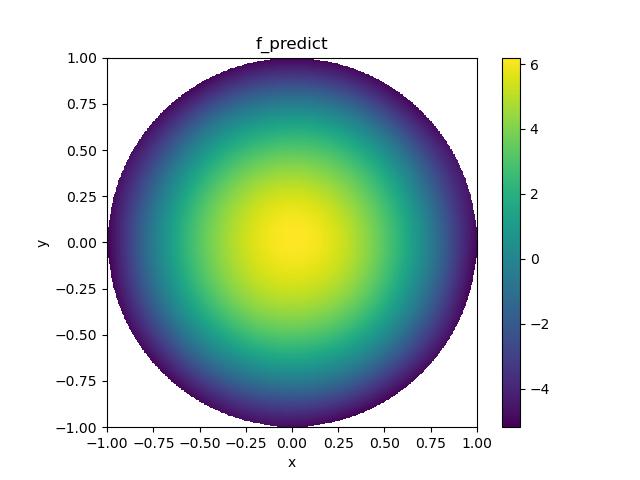}}
\subfloat[$\hat{f}$ with $\delta=10\%$.]
{\includegraphics[width=0.30\linewidth]{figures_2D_1.8/f_predict_2d_1.8_0.05.jpg}}
\\
\subfloat[$|\hat{f}-f^{*}|$ with $\delta=1\%$.]
{\includegraphics[width=0.30\linewidth]{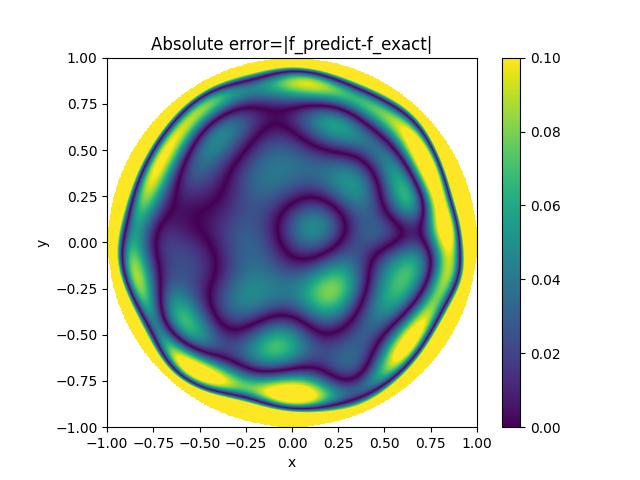}}
\subfloat[$|\hat{f}-f^{*}|$ with $\delta=5\%$.]
{\includegraphics[width=0.30\linewidth]{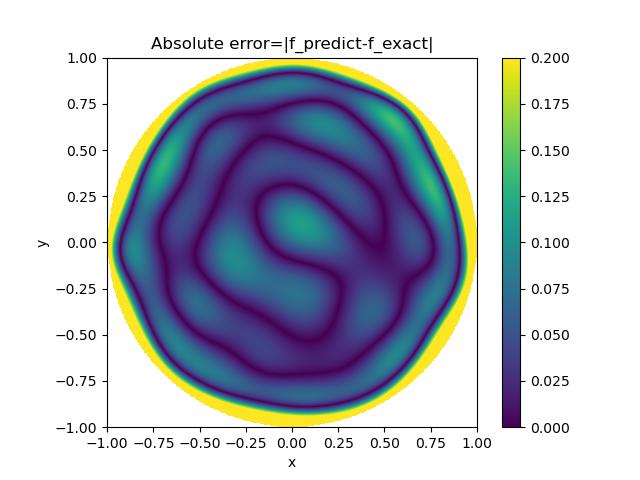}}
\subfloat[$|\hat{f}-f^{*}|$ with $\delta=10\%$.]
{\includegraphics[width=0.30\linewidth]{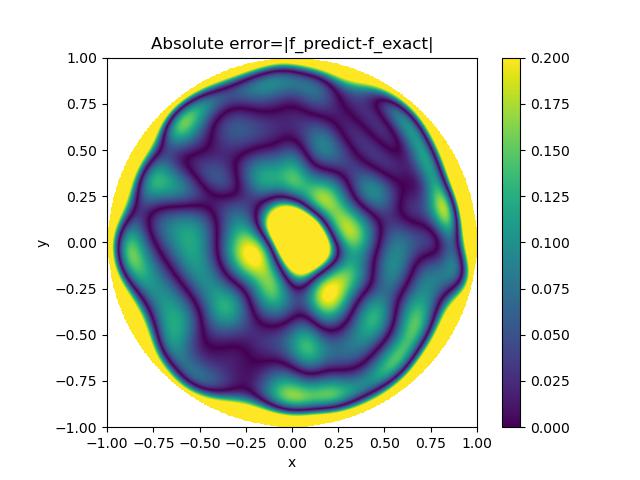}}
\caption{The reconstructions of $f$ (top) and the corresponding point-wise absolute error $|\hat{f}-f^{*}|$ (bottom) when $\alpha=1.8$}
\label{Figure 4}
\end{figure}

\subsection{Two-dimensional example} 

To begin with, we consider a 2-dimension problem. Since the noisy level of measurement data may influence the approximation accuracy of the neural network, we first report the relative errors for the reconstruction of $f$ and the absolute errors with respect to different fractional order $\alpha=0.5,1.2,1.5,1.8$ for different noisy levels: $\delta=0.01,0.05,0.1$.

According to Table \ref{Table 1}, the relative error of $\hat{f}$ decays when $\alpha$ decreases in [0,2], while the noise level has little influence on $Re_{f}$. We also plot the reconstruction results of $f(x)$ for four different choice of $\alpha$ in Figure \ref{Figure 1}-\ref{Figure 4}, in which the inversion result solved by the proposed method (first line) and the point-wise absolute error (second line) at different noisy level have been presented. As we can see, for different $\alpha$ the reconstructed $\hat{f}(x)$ is very close to the exact solution $f^*(x)$.

Next we take $\alpha=1.5$ and show the reconstructed solution for $f(x)$ and the corresponding absolute errors in higher dimension cases.

\subsection{High-dimensional example} 

\begin{figure}
\centering
\subfloat[$\hat{f}$ with $\delta=1\%$.]
{\includegraphics[width=0.30\linewidth]{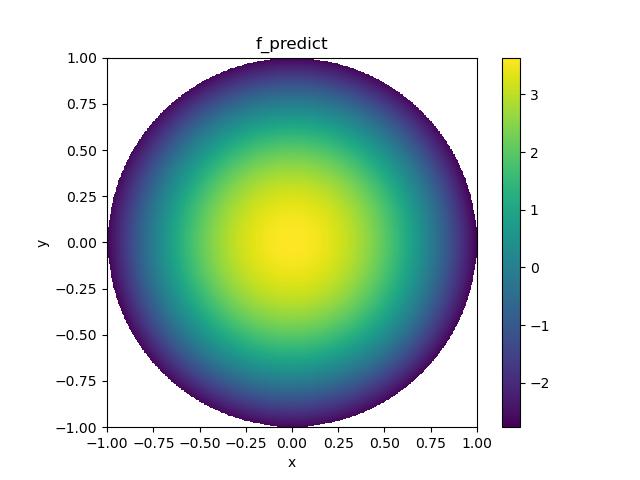}}
\subfloat[$\hat{f}$ with $\delta=5\%$.]
{\includegraphics[width=0.30\linewidth]{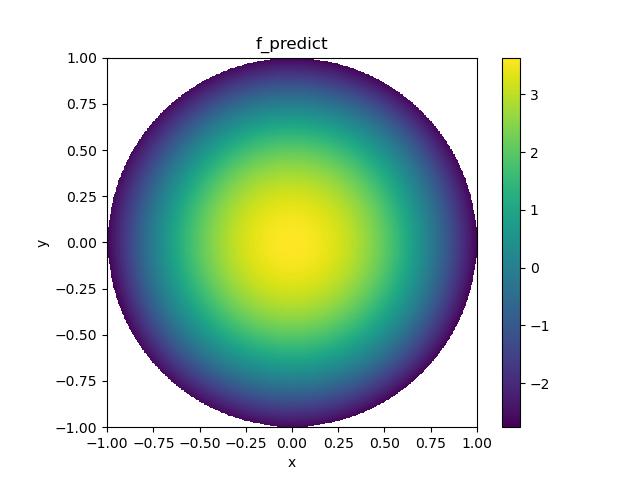}}
\subfloat[$\hat{f}$ with $\delta=10\%$.]
{\includegraphics[width=0.30\linewidth]{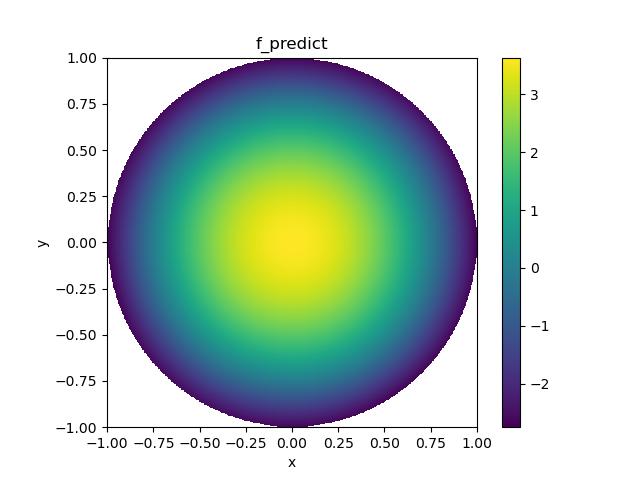}}
\\
\subfloat[$|\hat{f}-f^{*}|$ with $\delta=1\%$.]
{\includegraphics[width=0.30\linewidth]{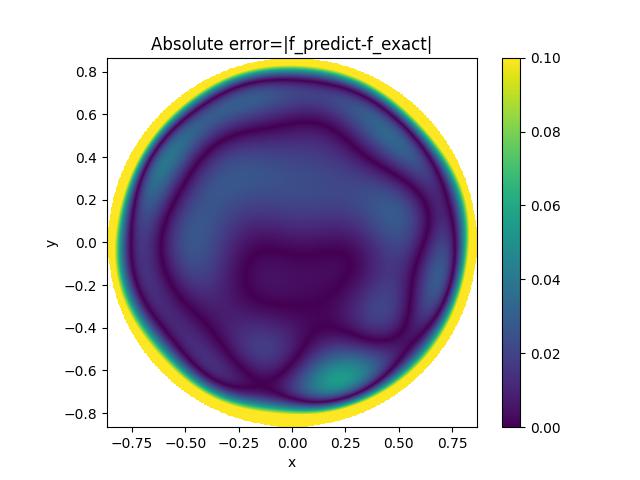}}
\subfloat[$|\hat{f}-f^{*}|$ with $\delta=5\%$.]
{\includegraphics[width=0.30\linewidth]{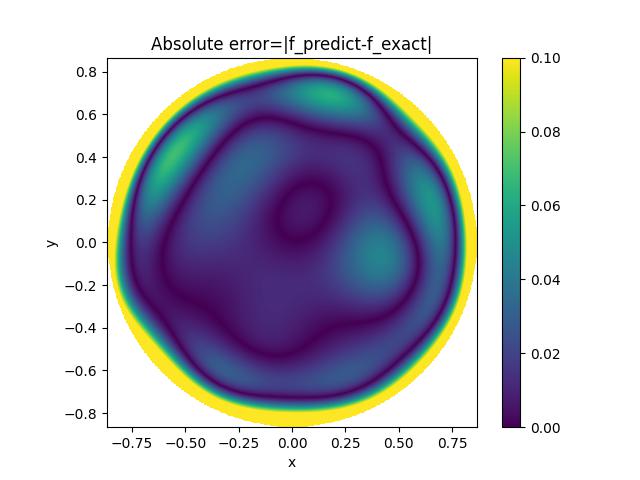}}
\subfloat[$|\hat{f}-f^{*}|$ with $\delta=10\%$.]
{\includegraphics[width=0.30\linewidth]{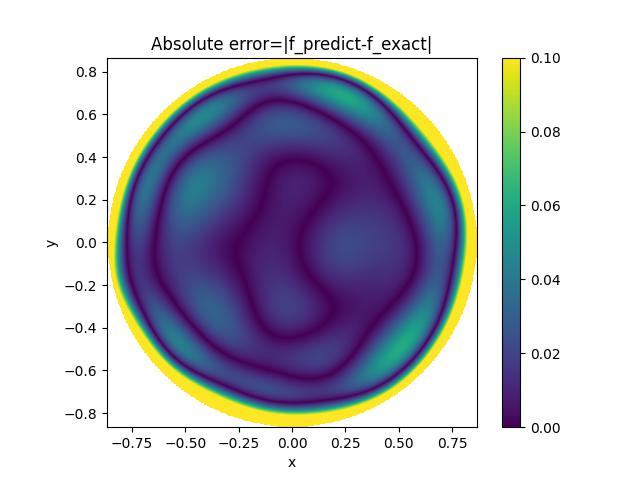}}
\caption{3D The reconstruction of $f$ and the corresponding point-wise absolute error $|\hat{f}-f^{*}|$ (bottom) at $x_3=0.5$ when $\alpha=1.5$}
\label{Figure 5}
\end{figure}

\begin{figure}
\centering
\subfloat[$\hat{f}$ with $\delta=1\%$.]
{\includegraphics[width=0.30\linewidth]{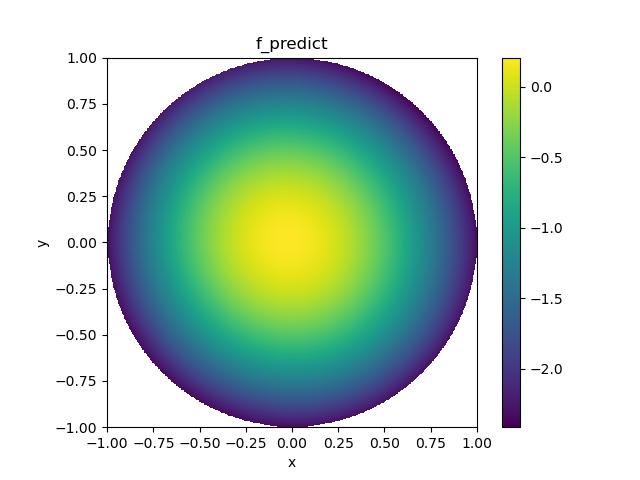}}
\subfloat[$\hat{f}$ with $\delta=5\%$.]
{\includegraphics[width=0.30\linewidth]{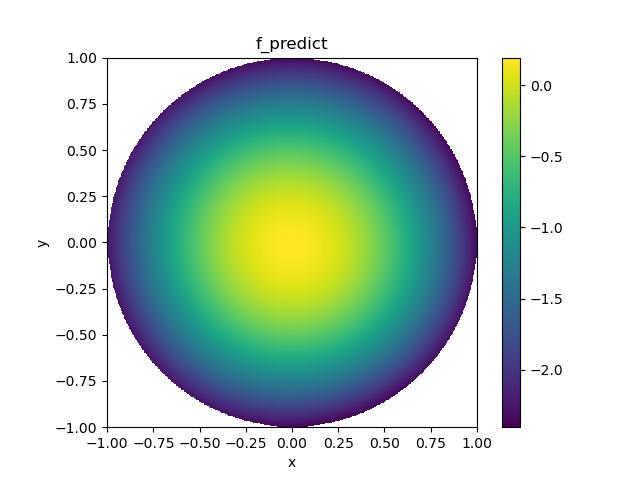}}
\subfloat[$\hat{f}$ with $\delta=10\%$.]
{\includegraphics[width=0.30\linewidth]{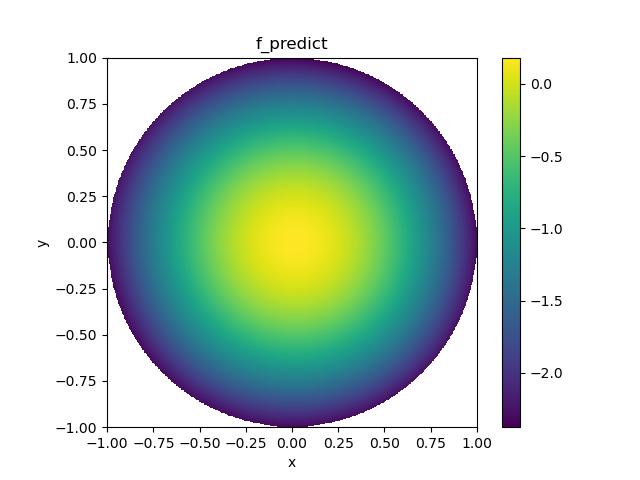}}
\\
\subfloat[$|\hat{f}-f^{*}|$ with $\delta=1\%$.]
{\includegraphics[width=0.30\linewidth]{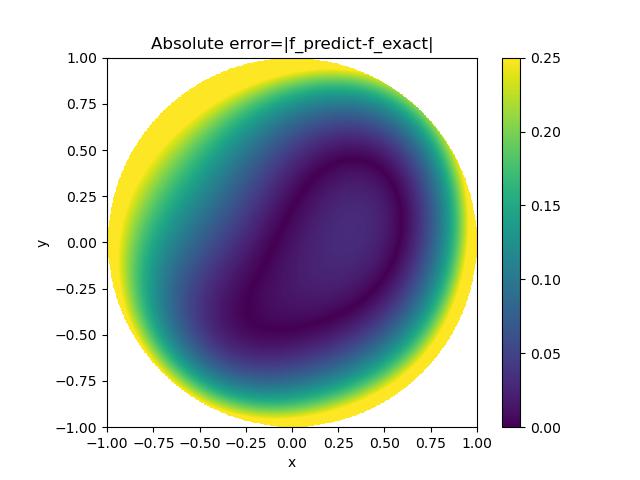}}
\subfloat[$|\hat{f}-f^{*}|$ with $\delta=5\%$.]
{\includegraphics[width=0.30\linewidth]{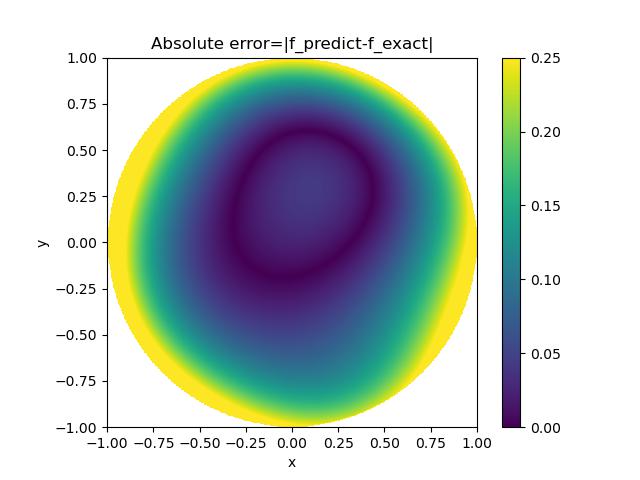}}
\subfloat[$|\hat{f}-f^{*}|$ with $\delta=10\%$.]
{\includegraphics[width=0.30\linewidth]{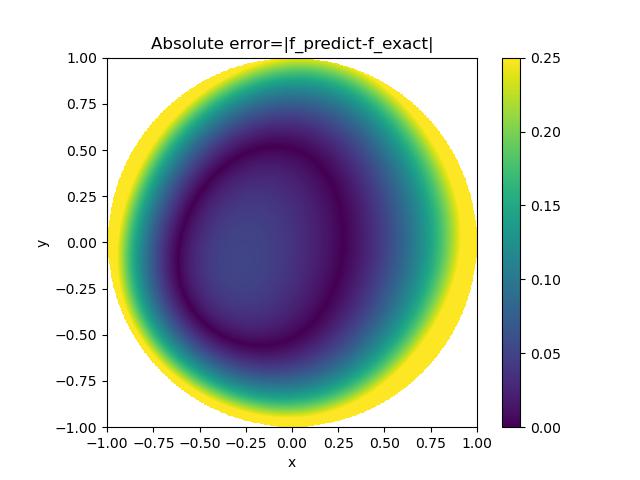}}
\caption{5D The reconstruction of $f$ and the corresponding point-wise absolute error $|\hat{f}-f^{*}|$ (bottom) when $\alpha=1.5$}
\label{Figure 6}
\end{figure}

Next, we consider several higher-dimensional problems with $\alpha=1.5$. Table \ref{Table2} demonstrates the results for $2D$, $3D$, $5D$ and $10D$ problems, from which we can find that the relative error of $\hat{f}$ increases with the dimension becoming larger, while all the relative error saturates around $10^{-2}$. Similar to the 2D problem, the noisy level $\delta$ still has little influence on $Re_{f}$ in 3D/5D/10D case. In Figure \ref{Figure 5} and Figure \ref{Figure 6}, we also plot the reconstruction results $\hat{f}(x)$ with its absolute errors for 3D and 5D problem respectively, which demonstrates the effectiveness of MC-fPINNs in dealing with such high dimensional problems.

\begin{table}
\centering
\begin{tabular}{cccc}
\toprule
&\multicolumn{3}{c}{$\delta$}   \\
\cmidrule(lr){2-4}
dimension & $1\%$ & $5\%$ & $10\%$  \\
\cmidrule{1-4}
2D & $2.25 \times 10^{-2}$ & $2.91 \times 10^{-2}$ & $3.17 \times 10^{-2}$ \\
3D & $2.81 \times 10^{-2}$ & $2.81 \times 10^{-2}$ &  $2.48 \times 10^{-2}$ \\
5D & $4.52 \times 10^{-2}$ & $4.82 \times 10^{-2}$ & $4.83 \times 10^{-2}$   \\
10D & $8.08 \times 10^{-2}$ & $8.45 \times 10^{-2}$ & $8.41 \times 10^{-2}$   \\
\bottomrule
\end{tabular}
\caption{The relative $L^{2}$-error of the reconstructed $\hat{f}$ in different dimension at different noise level when $\alpha=1.5$ }
\label{Table2}
\end{table}

\section{Conclusion}\label{conclusion}
In this paper, we study the application of MC-fPINNs in solving the inverse source problem of fractional Poisson equation. With a loss function consisting of the governing equation and measurement data, we utilize two neural networks simultaneously to approximate the solution $u^*$ and the forcing term $f^*$. Several numerical examples have been shown to illustrate the effectiveness of this method, including some higher-dimensinal problems. Furthermore, a rigorous error analysis of this method have been presented, providing us with a practical criterion for choosing parameters of neural networks. In our future work, we shall extend this method together with the error analysis to more complex fractional PDEs, such as the fractional advection-diffusion equation and the time-space fractional diffusion equation.

\section*{Acknowledgement}
This work is supported by the National Key Research and Development Program
of China (No.2020YFA0714200), by the National Nature Science Foundation of China (No. 12125103, No. 12071362, No. 12301558), and by the
Fundamental Research Funds for the Central Universities.

\bibliography{MCFP}
\bibliographystyle{plain}

\end{document}